# A Copula Statistic for Measuring Nonlinear Multivariate Dependence

Mohsen Ben Hassine[a], Lamine Mili[b*], Kiran Karra[c]

[a]Department of Computer Science, University of El Manar, Tunis, Tunisia. (E-mail: mbh851@yahoo.fr).

[b*]Corresponding author: Bradley Department of Electrical and Computer Engineering, Northern Virginia Center, Virginia Tech, Falls Church, VA 22043, USA (Tel: (703) 740 7610; E-mail: lmili@vt.edu).

[c]Bradley Department of Electrical and Computer Engineering, VTRC-A, Arlington, VA 22203, USA (E-mail: kiran.karra@vt.edu).

*Abstract—* **A new index based on empirical copulas, termed the Copula Statistic (*CoS*), is introduced for assessing the strength of multivariate dependence and for testing statistical independence. New properties of the copulas are proved. They allow us to define the CoS in terms of a relative distance function between the empirical copula, the Fréchet-Hoeffding bounds and the independence copula. Monte Carlo simulations reveal that for large sample sizes, the *CoS* is approximately normal. This property is utilised to develop a *CoS*-based statistical test of independence against various noisy functional dependencies. It is shown that this test exhibits higher statistical power than the Total Information Coefficient (*TICe*), the Distance Correlation (*dCor*), the Randomized Dependence Coefficient (*RDC*), and the Copula Correlation (*Ccor*) for monotonic and circular functional dependencies. Furthermore, the $R^2$-equitability of the *CoS* is investigated for estimating the strength of a collection of functional dependencies with additive Gaussian noise. Finally, the *CoS* is applied to a real stock market data set from which we infer that a bivariate analysis is insufficient to unveil multivariate dependencies and to two gene expression data sets of the Yeast and of the E. Coli, which allow us to demonstrate the good performance of the *CoS*.**

## I. INTRODUCTION

Measures of statistical dependence among random variables and signals are paramount in many scientific, engineering, signal processing, and machine learning applications. They allow us to find clusters of data points and signals, test for independence to make decisions, and explore causal relationships. The conventional measure of dependence is provided by the correlation coefficient, which was introduced in 1895 by Karl Pearson [1]. Since it relies on moments, it assumes statistical linear dependence. However, in biology, ecology and finance, to name a few, applications involving nonlinear multivariate dependence prevail. For such applications, the correlation coefficient is unreliable. Hence, several alternative metrics have been proposed over the last decades. Two popular rank-based metrics are the Spearman's $\rho_s$ [2] and Kendall's $\tau_k$ [3]. Modified versions of these statistics for monotonic dependence have been suggested in [4-6]. Although Spearman's rho, Kendall's tau, and related metrics work well in the monotonic case, they fail when non-monotonic dependencies arise.

Many proposals have been initiated by researchers in order to address this deficiency. Reshef *et al*. [7, 8] introduced the Maximal Information Coefficient (*MIC*) and the Total Information Coefficient (*TIC*), Lopes-Paz *et al*. [9] proposed the Randomized Dependence Coefficient (*RDC*), and Ding *et al*. [10, 11] put forth the Copula Correlation Coefficient (*Ccor*). Additionally, Székely *et al*. [12] proposed the distance correlation (*dCor*). These metrics are able to measure monotonic and nonmonotonic dependencies between random variables, but each has strengths and shortcomings. For instance, the *TICe*, *RDC* and *dCor* have high statistical power against independence among many different types of functional dependencies [9] than the other proposed statistics. Another measure of dependence based on the relative distance between copulas, termed the Target Dependence Coefficient (*TDC*), has been recently developed by Marti *et al.* [13]. It is defined as the ratio of the first Wasserstein distance between the independence and the empirical copula to the minimum of the first Wasserstein distances between the empirical copula and a predefined set of copulas.

Now, a question arises is what properties must a measure of association strength satisfy? Rényi [14] proposed a list of seven axioms, which include the axioms of symmetricity and functional invariance. Recently, Reshef *et al*. [7,8] proposed another important property, the $R^2$-equitability; it is defined as the ability of a statistic to equally measure the strength of a range of functional relationships between two random variables, *X* and *Y* satisfying $Y = f(X) + \varepsilon$, subject to the same level of additive noise $\varepsilon$ using the well-known squared Pearson correlation $R^2$. Later, Kinney and Atwal [15] argued that this definition is applicable to only some type of dependence of the noise, $\varepsilon$, on $f(X)$, is not invariant under invertible transformation of *X*, and relies on an unreliable measure of functional dependence, the $R^2$. Instead, they proposed another more general concept of equitability, the self-equitability, which makes use of the metric itself instead of the $R^2$ for assessing its equitability. They showed that the estimator of the *MIC* proposed in [16], termed the *MICe,* is not self-equitable when measuring bivariate dependence among random variables.

In signal processing, communications, and several engineering fields (e.g. power systems and finance), copulas (for an introduction see [17]) are gaining a great deal of attention to model nonlinear dependencies of stochastic signals (e.g., [18-23]). For instance, Krylov *et al.* [18] applied Gaussian, Student's t, and Archimedean copulas to supervised dual-polarization image classification in Synthetic Aperture Radar (SAR) while Sundaresan and Varshney [19] explored the use of copulas for estimating random signal sources from sensor observations. Recognizing the need of new copulas in many engineering applications, Zeng *et al.* [20, 21] derived the bivariate exponential, Rayleigh, Weibull, and Rician copulas, proved that the first three copulas are equivalent, and proposed a technique based on mutual information for selecting the best copula for a given dataset.

In this paper, we introduce a new index based on copulas, termed the *CoS*, for measuring the strength of nonlinear multivariate dependence and for testing for independence. To derive the *CoS*, new properties of the copulas are stated and proved. In particular, we show that a continuous bivariate copula simultaneously attains the upper and the lower Fréchet-Hoeffding bounds and the independence copula at a point of the unit square domain if and only if that point is a global optimum of the function that relates the two random variables associated with the copula. We also prove that a continuous bivariate copula attains either the upper or the lower Fréchet-Hoeffding bound over a domain if and only of there is a function that relates the two random variables associated with the copula, and that function is either non-decreasing or non-increasing over that domain, respectively. Based on these properties, we derive the *CoS* as a rank-based statistic that is expressed in terms of a relative distance function defined as the difference between the empirical copula and the product copula scaled by the difference between either the upper or the lower Fréchet-Hoeffding bounds and the product copula. We prove that the *CoS* ranges from zero to one and attains its lower and upper limit for the independence and the functional dependence case, respectively.

Monte Carlo simulations are carried out to estimate bias and standard deviation curves of the *CoS,* to assess its power when testing for independence, and to evaluate its equitability for functional dependence and its reliability for non-functional one. The simulations reveal that for large sample sizes, the *CoS* is approximately normal and approaches Pearson's $\rho_P$ for the Gaussian copula and Spearman's $\rho_S$ for the Gumbel, Clayton, Galambos, and BB6 copula. The *CoS* is shown to exhibit strong statistical power in various functional dependencies including linear, cubic, and fourth root dependence as compared to the *RDC, Ccor, MICe*, and the *dCor,* and outperforms the *MICe* for copula-induced dependence and for Ripley's forms. It also shows reasonable performance in non-monotonic functional dependencies including sinusoidal and quadratic, and exceptional performance in the circular functional dependence. Finally, the *CoS* is applied to real stock market returns; analysis shows that it performs as well as the *RDC*, *Ccor, MICe,* and the *dCor* in revealing bivariate dependencies. Additionally, the *CoS*'s unique ability to measure multivariate dependence is put forward by demonstrating that bivariate analysis is insufficient to fully analyze the financial returns data. The performance of the *CoS* is also compared to the *MICe, RDC, Ccor,* and the *dCor* using gene expression data sets of the *Yeast* and of the *E. Coli*. Receiver operating characteristic curves and F-scores show that the *CoS* performs well in all the tested cases.

The paper is organized as follows. Section II states Sklar's definition of a copula, recalls Fréchet's theorem, and proves two other theorems on copulas. Section III introduces a relative distance function and proves several of its properties. Section IV defines the *CoS* and provides an algorithm that implements it. Section V investigates the statistical properties of the CoS and describes a *CoS*-based statistical test of independence. Additionally, in this section, the $R^2$-equitability of the *CoS* is investigated for estimating the strength of a collection of functional dependencies with additive Gaussian noise. Section VI compares the performance of the *CoS* with the *dCor, RDC, Ccor,* and the *MICe* in measuring bivariate functional and non-functional dependencies between synthetic datasets, shows how the *CoS* can unveil multivariate dependence in real datasets of stock market returns, and demonstrates it good performance on gene expression networks.

## II. INTRODUCTION TO COPULAS

We first define a *d*-dimensional copula, then we recall Sklar's and Fréchet's theorems, and finally we prove two theorems that will allow us to define the *CoS*.

*Definition 1:* (Sklar [24]): A *d*-dimensional copula $C(u_1, ..., u_d)$ is a function that maps the unit hypercube, $\mathbf{I}^d = [0,1]^d$, to the unit interval, **I**, such that

a) $C(0, ..., 0, u_i, 0, ..., 0) = 0$ for $i = 1, ..., d$;

b) $C(1, ..., 1, u_i, 1, ..., 1) = u_i$ for $i = 1, ..., d$;

c) For all $(u_{11}, ..., u_{d1})$ and $(u_{12}, ..., u_{d2}) \in \mathbf{I}^d$ such that $u_{i1} \leq u_{i2}$ for $i = 1, ..., d$, we have the rectangular inequality,

$$\sum_{i_1=1}^{2} \cdots \sum_{i_d=1}^{2} (-1)^{i_1+\cdots+i_d} C(u_{1i_1}, ..., u_{di_d}) \geq 0. \tag{1}$$

The rectangle inequality in (1) ensures that the copula is a non-negative *d*-increasing function. Combined with the other two properties defined by Sklar, it satisfies all the properties of a joint probability distribution function of n random variables that are uniformly distributed over the unit hypercube [17]. Now the question that arises is the following: for a given joint probability distribution function, $H(x_1, ..., x_d)$ of *d* random variables, $X_1, ..., X_d$, and its associated marginal distribution functions, $F_1(x_1), ..., F_d(x_d)$, does there exist a copula that relates the joint to the marginal distribution functions? Is this copula unique? Both questions are addressed by Sklar's theorem [17, 24], which guarantees for any H(**x**) the existence of a copula C(.) defined as

$$H(x_1, \ldots, x_d) = C(F_1(x_1), \ldots, F_d(x_d)). \tag{2}$$

It also states that *C*(.) is unique if the random variables are continuous, and is uniquely determined over the product set of the ranges of the marginal distribution functions, otherwise. On the other hand, given a set of marginal distribution functions, $F_1(x_1), ..., F_d(x_d)$, and a *d*-dimensional copula, *C*(.), then *H*(.) given by (2) is the associated joint distribution function. It is apparent from (2) that a copula encompasses all the dependencies between the random variables, $X_1, ..., X_d$.

Another theorem [17] that we rely on to derive our *CoS* provides the Fréchet-Hoeffding lower and upper bounds of a copula $C(\mathbf{u})$ for any $\mathbf{u} \in \mathbf{I}^d$. This theorem states that

$$W(\mathbf{u}) \leq C(\mathbf{u}) \leq M(\mathbf{u}), \tag{3}$$

where

$$W(\mathbf{u}) = \text{Max}\{\sum_{i=1}^{d} u_i + 1 - d, \ 0\}, \tag{4}$$

and

$$M(\mathbf{u}) = P(U_1 \leq u_1, ..., U_n \leq u_d) = \text{Min}(u_1, ..., u_d). \tag{5}$$

Unlike $M(\mathbf{u})$, which is a copula for all $n \geq 2$, $W(\mathbf{u})$ is a copula for $n = 2$, but not for $n \geq 3$. However, as shown in [17], there exists a copula that is equal to $W(\mathbf{u})$ for $n \geq 3$. Another special case of interest is the one where the random variables are independent, yielding a product copula given by

$$\Pi(u_1,...,u_d) = u_1 u_2 ... u_d. \tag{6}$$

In the following, we first restrict attention to two-dimensional copulas to develop a statistical index, the *CoS*, in the bivariate dependence case and then we extend that index to multivariate dependence. To define the *CoS* of two continuous random variables $X$ and $Y$ with copula $C(u, v)$, we first provide three definitions of bivariate dependencies, from weaker to stronger versions, as introduced by Lehmann [25], then we state three theorems.

*Definition 2:* Two random variables, $X$ and $Y$, are said to be concordant (or discordant) if they tend to simultaneously take large (or small) values.

A more formal definition is as follows. Let $X$ and $Y$ be two random variables taking two pairs of values, $(x_i, y_i)$ and $(x_j, y_j)$. $X$ and $Y$ are said to be concordant if $(x_i - x_j)(y_i - y_j) > 0$; they are said to be discordant if the inequality is reversed.

*Definition 3:* Two random variables, $X$ and $Y$, defined on the domain $\mathfrak{D} = \text{Range}(X) \times \text{Range}(Y)$ are said to be Positively Quadrant Dependent (PQD) if

$$P(X \leq x, Y \leq y) \geq P(X \leq x) P(Y \leq y),$$

that is, $C(u,v) \geq \Pi(u,v)$ and Negative Quadrant Dependent (NQD) if

$$P(X \leq x, Y \leq y) < P(X \leq x) P(Y \leq y),$$

That is $C(u,v) \leq \Pi(u,v)$, for all $(x, y) \in \mathfrak{D}$.

*Definition 4:* Two random variables, $X$ and $Y$, are said to be comonotonic (respectively countermonotonic) if $Y = f(X)$ almost surely and $f(.)$ is an increasing (respectively a decreasing) function.

In short, we say that two random variables are monotonic if they are either comonotonic or countermonotonic.

*Theorem 1:* (Fréchet [17]: Let $X$ and $Y$ be two continuous random variables. Then, we have

a) $X$ and $Y$ are comonotonic if and only if the associated copula is equal to its Fréchet-Hoeffding upper bound, that is, $C(u, v) = M(u, v) = Min(u,v)$;

b) $X$ and $Y$ are countermonotonic if and only if the associated copula is equal to its Fréchet-Hoeffding lower bound, that is, $C(u, v) = W(u, v) = Max(u+ v -1, 0)$

c) $X$ and $Y$ are independent if and only if the associated copula is equal to the product copula, that is, $C(u,v) = \Pi(u,v)= uv$.

*Lemma 1:* Let $X$ and $Y$ be two continuous random variables with copula $C(F_1(x), F_2(y)) = H(x,y) = P(X \leq x, Y \leq y)$. Then we have

$$a) \ P(X \leq x, Y > y) = F_1(x) - C(F_1(x), F_2(y)); \tag{7}$$

$$b) \ P(X > x, Y \leq y) = F_2(y) - C(F_1(x), F_2(y)); \tag{8}$$

$$c) \ P(X > x, Y > y) = 1 - F_1(x) - F_2(y) + C(F_1(x), F_2(y)). \tag{9}$$

Proof: Let us partition the domain $\mathfrak{D} = \text{Range}(X) \times \text{Range}(Y)$ of the joint probability distribution function, $H(x,y)$, into four subsets, namely $\mathfrak{D}_1 = \{X \leq x, Y \leq y\}$, $\mathfrak{D}_2 = \{X \leq x, Y > y\}$, $\mathfrak{D}_3 = \{X > x, Y \leq y\}$ and $\mathfrak{D}_4 = \{X > x, Y > y\}$. Then we have

$$P(X \leq x, Y \leq y) + P(X \leq x, Y > y) + P(X > x, Y \leq y) + P(X > x, Y > y) = 1. \tag{10}$$

We also have

$$P(X \leq x, Y \leq y) + P(X \leq x, Y > y) = P(X \leq x), \tag{11}$$

which yields (7), and

$$P(X \leq x, Y \leq y) + P(X > x, Y \leq y) = P(Y \leq y), \tag{12}$$

which yields (8). Substituting the expressions of $P(X \leq x, Y > y)$ given by (7) and of $P(X > x, Y \leq y)$ given by (8) into (10), we get (9). ∎

In the following theorems and corollaries, we assume that $X$ and $Y$ are continuous random variables and related via a function $f(.)$, that is, $Y=f(X)$, where $f(\cdot)$ is continuous and differentiable over the range of $X$.

*Theorem 2:* Let $X$ and $Y$ be two continuous random variables such that $Y=f(X)$ almost surely, and let $C(u,v)$ be the copula value for the pair $(x,y)$. The function $f(.)$ has a global maximum at $(x_1, y_{max})$ with a copula value $C(u_1, v_1)$ or a global minimum at $(x_2, y_{min})$ with a copula value $C(u_2, v_2)$ if and only if we have

$$a) \ C(u_1, v_1) = M(u_1, v_1) = W(u_1, v_1) = \Pi(u_1, v_1) = u_1; \tag{13}$$

$$b) \ C(u_2, v_2) = M(u_2, v_2) = W(u_2, v_2) = \Pi(u_2, v_2) = 0. \tag{14}$$

Proof: a) Under the assumption that $Y=f(X)$, suppose that $(x_1, y_{max})$ is a global maximum of $f(.)$. Then, by definition we have $C(F_1(x_1), F_2(y_{max})) = P(X \leq x_1, Y \leq y_{max}) = P(X \leq x_1)$, implying that $C(u_1, 1) = u_1$. We also have $\text{Min}(u_1, v_1) = \text{Min}(u_1, 1) = u_1$ and $\text{Max}(u_1, v_1) = \text{Max}(u_1 + 1 - 1, 0) = u_1$, from which (13) follows. Let us prove the converse under the assumption that $Y=f(X)$. Suppose that there exists a pair $(u_1, v_1)$ such that $C(u_1, v_1) = M(u_1, v_1) = W(u_1, v_1) = \Pi(u_1, v_1) = u_1$. It follows that $v_1 = 1$, which implies that $C(u_1, 1) = u_1$ and $C(F_1(x_1), F_2(y_{max})) = P(X \leq x_1, Y \leq y_{max})$, that is, $(x_1, y_{max})$ is a global maximum of $f(.)$.

b) Suppose that $Y=f(X)$ and $(x_2, y_{min})$ is a global minimum. Then, by definition we have $C(F_1(x_2),F_2(y_{min})) = P(X \leq x_2, Y \leq y_{min}) = 0$, implying that $C(u_2,0) = 0$. We also have $W(u_2,v_2) = \min(u_2,0) = 0$, and $M(u_2, v_2) = \max(u_2 +0 -1,0) = 0$, from which (14) follows. Let us prove the converse under the assumption that $Y=f(X)$. Suppose that there exists a pair $(u_2, v_2)$ such that $C(u_2, v_2) = M(u_2, v_2) = W(u_2, v_2) = \Pi(u_2, v_2) = u_2 v_2 = 0$. It follows that either $u_2 = 0$, or $v_2 = 0$, or $u_2 = v_2 = 0$. Let us consider the first case where $u_2 = 0$. It follows that $C(0, v_2) = 0$, implying that $C(F_1(x_{2min}), F_2(y_2)) = P(X \leq x_{2min}, Y \leq y_2) = 0$. This means that $(x_{2min}, y_2)$ is a global minimum of $f(.)$. Let us consider the second case where $v_2 = 0$. It follows that $C(u_2, 0) = 0$, implying that $C(F_1(x_{2min}), F_2(y_2)) = P(X \leq x_2, Y \leq y_{2min}) = 0$. This means that $(x_2, y_{2min})$ is a global minimum of $f(.)$. Let us consider the third case where $u_2 = v_2 = 0$. It follows that $C(0, 0) = 0$, implying that $C(F_1(x_{2min}), F_2(y_{2\,min})) = P(X \leq x_{2\,min}, Y \leq y_{2min}) = 0$. This means that $(x_{2min}, y_{2\,min})$ is a global minimum of $f(.)$. ∎

*Corollary 1:* Let $X$ and $Y$ be two continuous random variables such that $Y = f(X)$, almost surely. If $f(.)$ is a periodic function, then (13) and (14) holds true at all the global maxima and global minima, respectively.

The proof of Corollary 1 directly follows from Theorem 2. This corollary is demonstrated in Fig. 1, which displays the graph of the projections on the $(u, C(u,v))$ plane of the empirical copula $C(u,v)$ associated with a pair $(X,Y)$, where $X$ is uniformly distributed over $[-1,1]$, and $Y = \sin(2\pi X)$. We observe that at each one of the four optima of the sine function, we have $C(u,v) = M(u,v) = W(u,v) = \Pi(u,v)$.

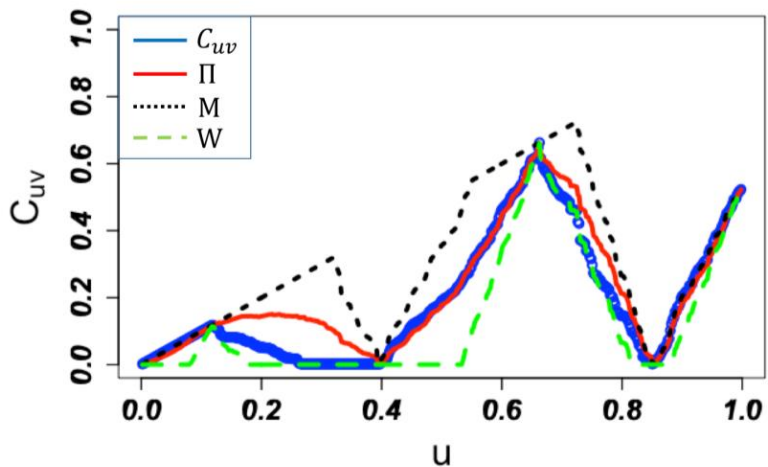

**Fig.1**. Graph (in blue dots) of the projections on the $(u, C(u,v))$ plane of the empirical copula $C(u,v)$ associated with a pair of random variables $(X, Y)$, where $X \sim U(-1,1)$ and $Y= \sin(2\pi\,X)$. The $u$ coordinates of the data points are equally spaced over the unity interval. Similar graphs are shown for the $M(u,v)$, $W(u,v)$ and $\Pi(u,v)$ copulas.

*Theorem 3:* Let *X* and *Y* be two continuous random variables such that *Y=f(X)*, almost surely where *f(.)* has a single optimum and let *C(u,v)* be the copula value for the pair *(x,y)*. We have $C(u,v) = M(u,v)$ if and only if $df(x)/dx \geq 0$ and $C(u,v)=W(u,v)$ otherwise.

Proof: Suppose that $Y = f(X)$ almost surely, where *f(.)* has a single optimum, which is necessarily a global one. Let us denote by $S_1$ and $S_2$ the non-increasing and the non-decreasing line segments of *f(.)*, respectively. Note that *f(.)* may have inflection points but may not have a line segment of constant value because otherwise *Y* will be a mixed random variable, violating the continuity assumption. Let *A* denote a point with coordinate *(x,y)* of the function *f(.)*. Consider the four subsets $\mathcal{D}_1 = \{X \leq x, Y \leq y\}$, $\mathcal{D}_2 = \{X \leq x, Y > y\}$, $\mathcal{D}_3 = \{X > x, Y \leq y\}$ and $\mathcal{D}_4 = \{X > x, Y > y\}$. Suppose that *A* is a point of $S_1$. As shown in Fig. 1(a), either $\mathcal{D}_1 \cap S_1 = \{A\}$ or $\mathcal{D}_4 \cap S_1 = \emptyset$ depending upon whether *f(.)* has a global minimum or a global maximum point, respectively. In the former case, we have $P(X \leq x, Y \leq y) = 0$, implying that $C(u,v) = 0$, while in the latter case, we have $P(X > x, Y > y) = 0$, implying from (9) that $C(u, v) = u + v - 1 \geq 0$. Combining both cases, it follows that for all $(x, y) \in S_1$, $C(u,v) = \text{Max}(u + v - 1, 0)$.

Now, suppose that *A* is a point of $S_2$. As shown in Fig. 2(b), either $\mathcal{D}_2 \cap S_2 = \{A\}$ or $\mathcal{D}_3 \cap S_2 = \emptyset$ depending upon whether *f(.)* has a global maximum or a global minimum point, respectively. In the former case, we have $P(X \leq x, Y > y) = 0$, implying from (7) that $C(u,v) = u$ while in the latter case, we have $P(X > x, Y \leq y) = 0$, implying from (8) that $C(u, v) = v$. Combining both cases, it follows from (3) that for all $(x,y) \in S_2$, $C(u,v) = \min(u,v)$. ∎

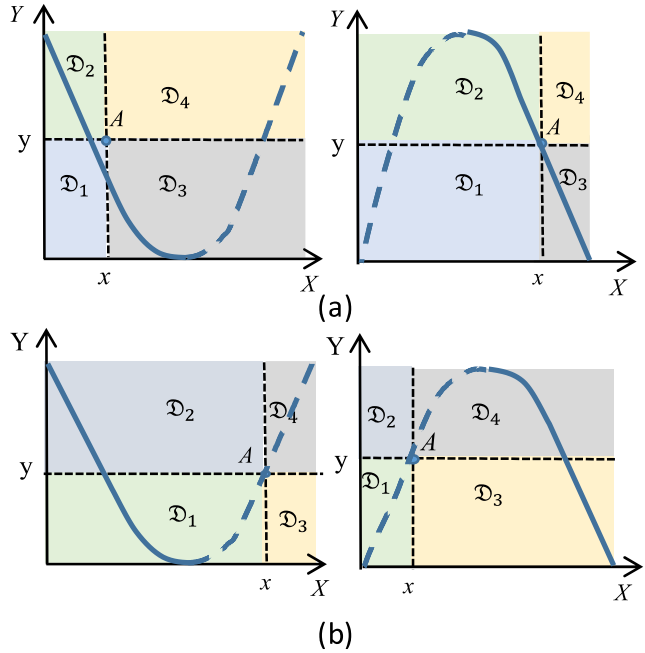

**Fig. 2**. Graphs of a function $Y = f(X)$ having a single optimum. A point A with coordinate *(x,y)* is located either on the non-increasing part, $S_1$, shown as a solid line in (a) or on the non-decreasing part, $S_2$, shown as a dashed line in (b) of the function $f(.)$. Four domains, $\mathcal{D}_1, \ldots, \mathcal{D}_4$, are delineated by the vertical and horizontal lines at position $X = x$ and $Y = y$, respectively.

Theorem 3 is illustrated in Fig. 3. This figure displays the graph of the projections on the (u, C(u,v)) plane of C(u,v) associated with a pair of random variables, (X,Y), where X follows U(-5, 5) and $Y = f(X) = (X-1)^2$. We observe that C(u,v) = W(u,v) for $0 \leq u \leq 0.6$, for which $f'(x) \leq 0$ and C(u,v) = M(u,v) for $0.6 \leq u \leq 1$, for which $f'(x) \geq 0$.

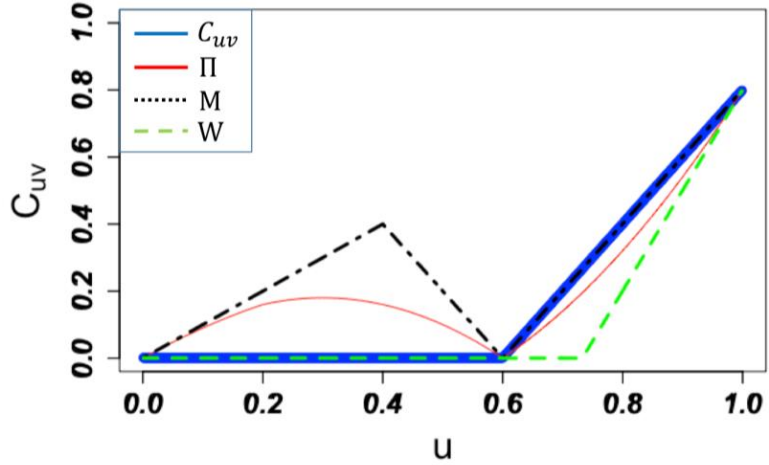

**Fig.3**. Graph (blue circles) of the projections on the (u, C(u,v)) plane of C(u,v) associated with $X \sim U(-5,5)$ and $Y = f(X) = (X-1)^2$. The u coordinates of the data points are equally spaced. The minimum of the function f(.) is associated with u = 0.6 and C(u,v) = 0. Similar graphs are shown for M(u,v) (dotted black), W(u,v) (dashed green), and Π(u,v) (solid red). We have C(u,v) = W(u,v) for $0 \leq u \leq 0.6$, which corresponds to $f'(x) \leq 0$, and C(u,v) = M(u,v) for $0.6 \leq u \leq 1$, which corresponds to $f'(x) \geq 0$.

## III. THE RELATIVE DISTANCE FUNCTION

We define a metric of proximity of the copula to the upper or the lower bounds with respect to the Π copula and investigate its properties.

*Definition 5:* The relative distance function, $\lambda(C(u,v))$: $[0,1] \to [0,1]$, is defined as

a) $\lambda(C(u,v)) = (C(u,v) - uv)/(\text{Min}(u,v) - uv)$ if $C(u,v) \geq uv$;

b) $\lambda(C(u,v)) = (C(u,v) - uv)/(\text{Max}(u+v-1,0) - uv)$ if $C(u,v) < uv$.

In other words, $\lambda(C(u,v))$ is equal to the ratio of the difference between $C(u,v)$ and $\Pi(u,v)$ to the difference between $M(u,v)$ (respectively $W(u,v)$) and $\Pi(u,v)$ if $X$ and $Y$ are PQD (respectively NQD). This is illustrated in Fig. 4. Note that we have $\lambda(C(u,v)) = 1$ if $C(u,v) = M(u,v)$ or $C(u,v) = W(u,v)$ and from (3), we have $W(u,v) \leq \Pi(u,v) \leq M(u,v)$.

*Theorem 4:* $\lambda(C(u,v))$ satisfies the following properties:

a) $0 \leq \lambda(C(u,v)) \leq 1$ for all $(u,v) \in \mathbf{I}^2$;

b) $\lambda(C(u,v)) = 0$ for all $(u,v) \in \mathbf{I}^2$ if and only if $C(u,v) = uv$;

c) If $Y = f(X)$ almost surely, where $f(.)$ is monotonic, then $\lambda(C(u,v)) = 1$ for all $(u,v) \in \mathbf{I}^2$;

d) If $Y = f(X)$ almost surely, then $\lambda(C(u,v)) = 1$ at the global optimal points of $f(.)$.

Proof: Property a) follows from Definition 5 and (3) while properties b), c) and d) follow from Definition 5 and Theorem 1 and 2. ∎

*Corollary 2:* If $Y = f(X)$ almost surely, where $f(\cdot)$ has a single optimum, then $\lambda(C(u,v)) = 1$ for all $(u,v) \in \mathbf{I}^2$.

Proof: it directly follows from Theorem 3 and Definition 5. ∎

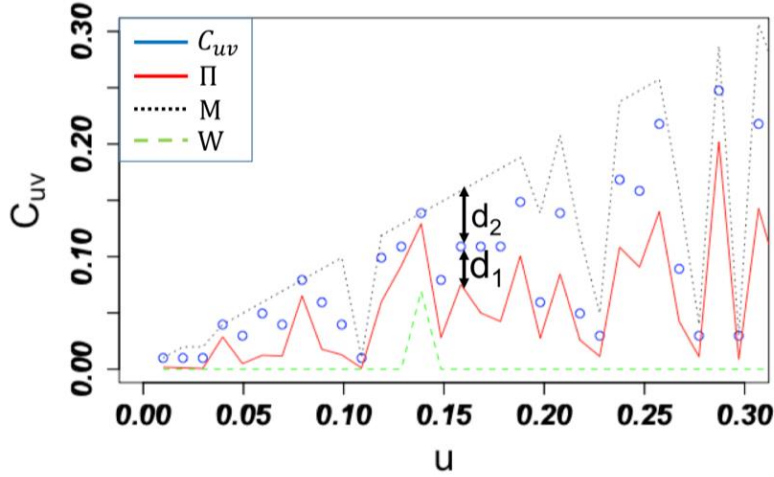

**Fig. 4.** Graph (blue circles) of the projections on the $(u, C(u,v))$ plane drawn from the Gaussian copula $C(u,v)$ with $\rho_P = 0.5$. Similar graphs are shown for $M(u,v)$ (dotted black), $W(u,v)$ (dashed green), and $\Pi(u,v)$ (solid red). The empirical relative distance function is given by $\lambda(C(u,v)) = d_1/d_2$, where $d_1$ is the distance from $C(u,v)$ to $\Pi(u,v)$ and $d_2$ is the distance from $M(u,v)$ to $\Pi(u,v)$.

Now, the question that arises is the following: Is $\lambda(C(u,v)) = 1$ for all $(u,v) \in \mathbf{I}^2$ when there is a functional dependence with multiple optima, be they global or local? The answer is given by the following two theorems.

*Theorem 5:* If $Y = f(X)$ almost surely where $f(\cdot)$ has at least two global maxima or two global minima and no local optima on the domain $\mathcal{D} = \text{Range}(X) \times \text{Range}(Y)$, then there exists a non-empty interval of $X$ for which $\lambda(C(u,v)) < 1$.

Proof: Suppose that $Y = f(X)$ almost surely, where $f(.)$ has at least two global maxima and no local optima. As depicted in Fig. 5(a), let $B$ and $C$ be two global maximum points of $f(.)$ with coordinates $(x_B, y_{\max})$ and $(x_C, y_{\max})$, respectively. This means that there exists $\Delta x > 0$ such that $f(x_B \pm \Delta x) < y_{\max}$ and $f(x_C \pm \Delta x) < y_{\max}$. Consider a point $A$ with coordinate $(x_A, y_A)$ such that $x_B < x_A < x_B + \Delta x$, $f(x_B - \Delta x) < y_A < y_{\max}$ and $f(x_C - \Delta x) < y_A < y_{\max}$. Let us denote by $S_B$ and $S_C$ the line segments of $f(.)$ defined over the intervals $[f(x_B - \Delta x), y_{\max}]$ and $[f(x_C - \Delta x), y_{\max}]$, respectively, which are shown as solid lines in Fig. 5(a). Let us partition the domain $\mathcal{D}$ into four subsets, $\mathcal{D}_1 = \{X \leq x_A, Y \leq y_A\}$, $\mathcal{D}_2 = \{X \leq x_A, Y > y_A\}$, $\mathcal{D}_3 = \{X > x_A, Y \leq y_A\}$ and $\mathcal{D}_4 = \{X > x_A, Y > y_A\}$. As observed in Fig. 3(a), we have $\mathcal{D}_1 \cap S_B\setminus\{A\} \neq \emptyset$, $\mathcal{D}_2 \cap S_B \neq \emptyset$, $\mathcal{D}_3 \cap S_C \neq \emptyset$, and $\mathcal{D}_4 \cap S_C \neq \emptyset$, yielding $\lambda(C(u,v)) < 1$. A similar proof can be developed for the case where $f(.)$ has at least two global minima and no local optima. ∎

Next, we prove a theorem that states that $\lambda(C(u,v))$ may be smaller than one at a local optimum of $f(.)$. Therefore, when developing the algorithm that implements the *CoS,* we must include a procedure that identifies all local optima of $f(.)$ and that sets the *CoS* equal to one at these points. This is achieved in Step 7 of the algorithm described in Section IV.C.

*Theorem 6:* If $Y = f(X)$ almost surely, where $f(\cdot)$ has a local optimum, then $\lambda(C(u,v)) \leq 1$ at that point.

Proof: Suppose that $Y = f(X)$ almost surely, where $f(.)$ has a local minimum point, say point A of coordinates $(x_A, y_A)$ as shown in Fig. 5(b). This means that there exists $\Delta x > 0$ such that $f(x_A \pm \Delta x) > y_A$. As depicted in Fig. 3(b), let $S_{A1}$ and $S_{A2}$ denote the line segments of $f(.)$ defined over $x_A - \Delta x$ and $x_A + \Delta x$, respectively. Let us consider the four domains, $\mathcal{D}_1 = \{X \leq x_A, Y \leq y_A\}$, $\mathcal{D}_2 = \{X \leq x_A, Y > y_A\}$, $\mathcal{D}_3 = \{X > x_A, Y \leq y_A\}$ and $\mathcal{D}_4 = \{X > x_A, Y > y_A\}$. As observed in Fig. 3(b), we have $\mathcal{D}_2 \cap S_{A1} \neq \emptyset$ and $\mathcal{D}_4 \cap S_{A2} \neq \emptyset$. Now, because A is by hypothesis a local minimum point, there exist line segments of $f(.)$ denoted by $S$ such that $f(y) < y_A$. Consequently, we have one of the following three cases: either $\mathcal{D}_1 \cap S\setminus\{A\} \neq \emptyset$ and $\mathcal{D}_3 \cap S \neq \emptyset$ as depicted in Fig. 3(b), or $\mathcal{D}_1 \cap S\setminus\{A\} \neq \emptyset$ and $\mathcal{D}_3 \cap S = \emptyset$, or $\mathcal{D}_1 \cap S\setminus\{A\} = \emptyset$ and $\mathcal{D}_3 \cap S \neq \emptyset$. In the first case, $\lambda(C(u,v)) < 1$ while in the last two cases, $\lambda(C(u,v)) = 1$. A similar proof can be developed for $f(.)$ with a local maximum point. ∎

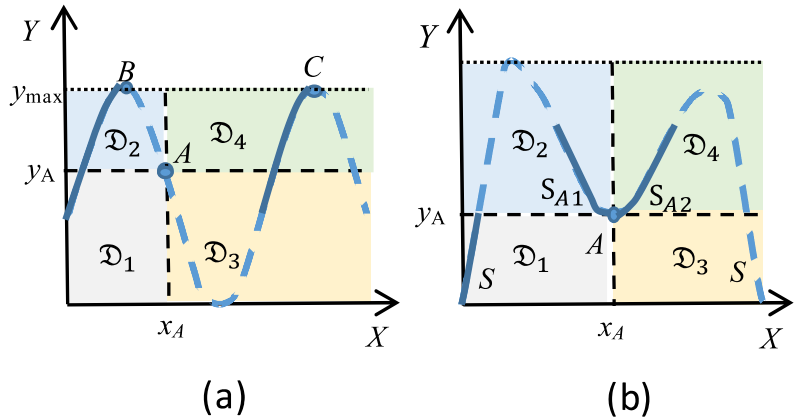

**Fig 5.** (a) The graph of a function $Y = f(X)$ having two global maximum points denoted by B and C, and one global minimum point, with two solid line segments denoted by $S_B$ and $S_C$. (b) The graph of a function $Y = f(X)$ having one local minimum point denoted by A, with line segments denoted by $S_{A1}$, $S_{A2}$, and $S$. Four domains, $\mathcal{D}_1, ..., \mathcal{D}_4$, are delineated by the vertical and horizontal lines at position $X = x_A$ and $Y = y_A$, respectively.

## IV. THE COPULA STATISTIC

We first define the empirical copula, then we introduce the copula statistic, and finally we provide an algorithm that implements it. One possible definition for the *CoS* is the mean of $\lambda(C(u,v))$ over $\mathbf{I}^2$, that is, $CoS(X,Y) = E[\lambda(C(u,v))]$. However, according to Theorems 5 and 6, $CoS \leq 1$ for functional dependence with multiple optima, which is not a desirable property. This prompts us to propose a better definition of the *CoS* based on the empirical copula as explained next.

### A. The Empirical Copula

Let $\{(x_i, y_i), i=1,\ldots, n, n \geq 2\}$ be a 2-dimensional data set of size *n* drawn from a continuous bivariate joint distribution function, $H(x, y)$. Let $R_{xi}$ and $R_{yi}$ be the rank of $x_i$ and of $y_i$, respectively. Deheuvels [26] defines the associated empirical copula as

$$C_n(u,v) = \frac{1}{n}\sum_{i=1}^{n} \mathbf{1}(u_i = \frac{R_{xi}}{n} \leq u, v_i = \frac{R_{yi}}{n} \leq v), \tag{15}$$

and shows its consistency. Here $\mathbf{1}(u_i \leq u, v_i \leq v)$ denotes the indicator function, which is equal to 0 or 1 if its argument is false or true, respectively. The empirical relative distance, $\lambda(C_n(u,v))$, satisfies Definition 4 by replacing $C(u,v)$ with the empirical copula given by (15).

### B. Defining the Copula Statistic for Bivariate Dependence

Let *X* and *Y* be two continuous random variables with a copula $C(u,v)$. Consider the ordered sequence, $x_{(1)} \leq \ldots \leq x_{(n)}$, of *n* realizations of *X*. This sequence yields $u_{(1)} \leq \ldots \leq u_{(n)}$ since $u_i = R_{xi}/n$ as given by (15). Let $\mathfrak{D}$ be the set of *m* contiguous domains $\{\mathfrak{D}_i, i = 1, \ldots, m\}$, where each $\mathfrak{D}_i$ is a *u*-interval associated with a non-decreasing or non-increasing sequence of $C_n(u_{(i)}, v_j)$, $i = 1, \ldots, n$. These domains form a partition of $\mathfrak{D}$, that is, $\mathfrak{D} = \cup_{i=1}^{m} \mathfrak{D}_i$ and $\mathfrak{D}_i \cap \mathfrak{D}_j = \emptyset$ for $i \neq j$. Let $C_i^{min}$ and $C_i^{max}$ respectively denote the smallest and the largest value of $C_n(u,v)$ on the domain $\mathfrak{D}_i$. Let $\gamma_i$ be defined as

$$\gamma_i = \begin{cases} 1 & \text{at a local optimum of } Y = f(X) \text{ on } \mathfrak{D}_i, \\ \frac{\lambda(C_i^{min}) + \lambda(C_i^{max})}{2}, & \text{otherwise.} \end{cases} \tag{16}$$

Note that the condition stated in (16) ensures that $\gamma_i = 1$ at a local optimum in the functional dependence case. We are now in a position to define the *CoS*.

*Definition 6:* Let $n_i$ denote the number of data points in the *i*-th domain $\mathfrak{D}_i$, $i = 1,\ldots, m$, while letting a boundary point belong to two contiguous domains, $\mathfrak{D}_i$ and $\mathfrak{D}_{i+1}$. Then, the copula statistic is defined as

$$CoS(X,Y) = \frac{1}{n+m-1}\sum_{i=1}^{m} n_i \gamma_i. \tag{17}$$

Note that we have $\sum_{i=1}^{m} n_i = n + m - 1$, yielding $CoS = 1$ if $\gamma_i = 1$ for $i = 1,\ldots, m$.

*Corollary 3:* The *CoS* of two random variables, *X* and *Y,* has the following asymptotic properties:

a) $0 \leq CoS(X,Y) \leq 1$;

b) $CoS(X,Y) = 0$ if and only if *X* and *Y* are independent;

c) If $Y = f(X)$ almost surely, then $CoS(X,Y) = 1$.

Proof: Properties a) and b) follow from Theorem 4 and the definitions given by (16) and (17). Property c) follows from Theorems 5, 6 and the definitions given by (16) and (17). ∎

Corollary 3c) states that $CoS(X,Y) = 1$ asymptotically for all types of functional dependence, which is a desirable property. What about the finite-sample properties of the *CoS* ? They are investigated in Section V. But first, let us describe an algorithm that implements the *CoS*.

*C. Algorithmic Implementation of the Copula Statistic*

Given a two-dimensional data sample of size *n*, $\{(x_j, y_j), j=1,..., n, n \geq 2\}$, the algorithm that calculates the *CoS* consists of the following steps:

1. Calculate $u_j$, $v_j$ and $C_n(u,v)$ as follows:

    a. $u_j = \frac{1}{n}\sum_{j=1}^{n} \mathbf{1}\{k \neq j : x_k \leq x_j\}$;

    b. $v_j = \frac{1}{n}\sum_{j=1}^{n} \mathbf{1}\{k \neq j : y_k \leq y_j\}$;

    c. $C_n(u,v) = \frac{1}{n}\sum_{j=1}^{n} \mathbf{1}\{u_j \leq u, v_j \leq v\}$;

2. Order the $x_j$'s to get $x_{(1)} \leq ... \leq x_{(n)}$, which results in $u_{(1)} \leq ... \leq u_{(n)}$ since $u_j = R_{x_j}/n$, where $R_{x_j}$ is the rank of $x_j$;

3. Determine the domains $\mathfrak{D}_i$, $i = 1, ... , m$, where each $\mathfrak{D}_i$ is a *u*-interval associated with a non-decreasing or non-increasing sequence of $C_n(u_{(j)},v_p), j = 1, ... , n$.

4. Determine the smallest and the largest value of $C_n(u,v)$, denoted by $C_i^{min}$ and $C_i^{max}$, and find the associated $u_i^{min}$ and $u_i^{max}$ for each domain $\mathfrak{D}_i$, $i = 1, ... , m$.

5. Calculate $\lambda(C_i^{min})$ and $\lambda(C_i^{max})$;

6. If $\lambda(C_i^{min})$ and $\lambda(C_i^{max})$ are equal to one, go to step 8;

7. Calculate the absolute difference between the three consecutive values of $C_n(u_{(i)},v_j)$ centered at $u_i^{min}$ (respectively at $u_i^{max}$) and decide that the central point is a local optimum if (i) both absolute differences are smaller than or equal to $1/n$ and (ii) there are more than four points within the two adjacent domains, $\mathfrak{D}_i$ and $\mathfrak{D}_{i+1}$;

8. Calculate $\gamma_i$ given by (16);

9. Repeat Steps 2 through 7 for all the *m* domains, $\mathfrak{D}_i$, $i = 1, ... , m$;

10. Calculate the *CoS* given by (17).

Note that Step 1 is the computation of the empirical copula as defined by Deheuvels [26]; steps 2-10 then utilize the empirical copula to compute the *CoS*. Step 7 checks whether a boundary point of a domain $\mathfrak{D}_i$ is a local optimum of $Y = f(X)$ and ensures that $\gamma_i = 1$ if that is the case. This rule is based on the following conjecture: $C_n(u_{(j)}, v_p)$ reaches a maximum (respectively a minimum) at a pair $(u_{(j)}, v_p)$ where $f(.)$ has a local maximum (respectively a local minimum). This conjecture stems from the extensive simulations that we carried out. The simulations also reveal that the variability of $C_n(u_{(j)}, v_p)$ vanishes when $X$ and $Y$ are functionally dependent, hence the test (ii) in Step 7. This is illustrated in Fig. 6 with a 4$^{th}$-order polynomial dependence having two global minima and one local maximum. As observed, at the global optimum points of $f(.)$ we have $C_n(u, v) = \pi(u, v) = M(u, v) = W(u, v)$, yielding $\lambda(C(u,v)) = 1$, while at the local maximum point of $f(.)$ we have $C_n(u, v) = \pi(u, v) \neq M(u, v) \neq W(u, v)$, yielding $\lambda(C(u,v)) < 1$.

### D. Defining the Multivariate CoS

Measure of multivariate dependence is receiving a growing attention in the literature [27-30]. Joe [29] is the first to extend Kendall's $\tau$ and Spearman $\rho_S$ to multivariate dependence. Following this development, Jouini and Clemen [30] propose a general expression for multivariate Kendall's $\tau$ based on the *d*-dimensional copula, $C(u_1, ..., u_d)$, which is defined as

$$\tau_n = \frac{1}{2^{d-1}-1} \left[ 2^d \int_{I^d} C(u_1, ..., u_d) \, dC(u_1, ..., u_d) - 1 \right]. \tag{18}$$

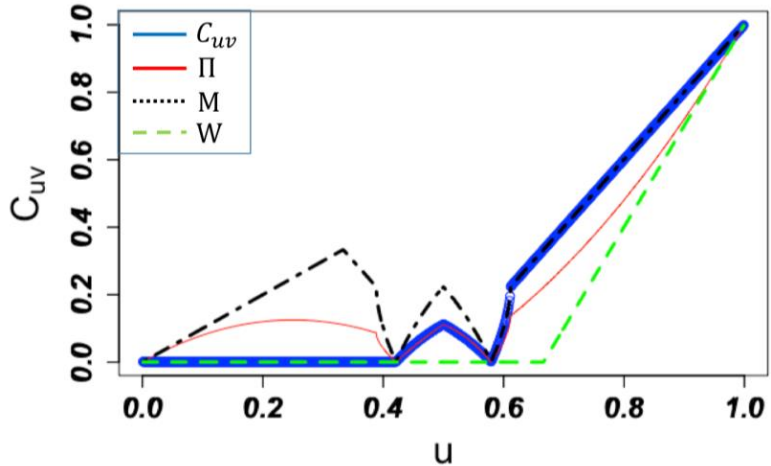

**Fig. 6.** Graph (blue circles) of the projections on the $(u, C(u,v))$ plane of $C(u,v)$ associated with $X \sim U(-5,5)$ and $Y = f(X) = (X^2 - 0.25)(X^2 - 1)$, which has two global minima at $x = \pm\sqrt{0.625}$ and one local maximum at $x = 0$. Similar graphs are displayed for $M(u,v)$ (dotted black), $W(u,v)$ (dashed green), and $\Pi(u,v)$ (solid red). The local optimum of $f(X)$ is associated with local optima of $C(u,v)$ and $\Pi(u,v)$ of equal magnitude shown at $u = 0.5$ on the graph.

While a multivariate version of the *MIC* has not been proposed yet, a multivariate *CoS* can be straightforwardly defined by extending the relative distance given by Definition 5 to the *d*-dimensional copula and the algorithm that implements the *CoS* given in Section IV.C to the *d*-dimensional empirical copula, which is expressed as

$$C_n(u_1, u_2, \ldots, u_d) = \frac{1}{n} \sum_{j=1}^{n} \mathbf{1}\{u_{1j} \leq u_1, \ldots, u_{dj} \leq u_d\};$$

where

$$u_{kj} = \frac{1}{n} \sum_{j=1}^{n} \mathbf{1}\{k \neq j : x_k \leq x_{kj}\}.$$

Unlike the *RDC* and *dCor*, which compute a metric of bivariate dependence between random vectors, the *CoS* can compute a metric of dependence between multiple serially dependent stochastic signals simultaneously. Although the number of copulas grows combinatorically as the dimensionality of the dataset increases, the *CoS* allows us to conduct multivariate analysis whenever possible, which may reveal further dependencies not uncovered by bivariate analysis. This fact will be demonstrated with financial returns data in Section IX.C.

### E. Computational Complexity of Calculating the CoS

The computational complexity of the algorithm described in Section IV.C for calculating the *CoS* is on the order of $O(d\, n\, log(n) + n^2 + n)$, where *d* is the dimensionality of the data being analyzed and *n* is the number of samples available to process. Specifically, $O(d\, n\, log\, n)$ is the run time complexity of the sorting operation involved in the computation of $u_{kj}$ for each dimension, $O(n^2)$ is the run time complexity of computing the empirical copula function, and $O(n)$ is the run time complexity of Step 2 to Step 10 of the algorithm. It is interesting to note that the run time complexity of the algorithm scales linearly with the dimension *d*, which allows us to compute the multivariate *CoS* in high dimensions (e.g., *d* > *10*).

## V. STATISTICAL PROPERTIES AND EQUITABILITY OF THE *CoS* FOR BIVARIATE DEPENDENCE

We analyze the finite-sample bias of the *CoS* for the independence case, then we develop a statistical test of bivariate independence, and finally we assess the $R^2$-equitability of the *CoS*.

### A. Statistical Analysis of the CoS

#### 1) Finite-Sample Bias of the CoS

Table I displays the sample means and the sample standard deviations of the *CoS* for independent random samples of increasing size generated from three monotonic copulas, namely Gauss(0), Gumbel(1), and Clayton(0), where a copula parameter value is indicated in brackets. As observed, the *CoS* has a bias for small to medium sample sizes. Interestingly, very close bias curves whose differences do not exceed 1% have been estimated from random samples drawn from a collection of 23 copulas using the copula package available on the CRAN repository website [31]. Fig. 7(a) shows a bias curve given by $CoS = 8.05\, n^{-0.74}$, fitted to 19 mean bias values for Gauss(0) using the least-squares method applied to a power model. It is observed that the *CoS* bias becomes negligible for a sample size larger than 500. Fig. 7(b) shows values taken by the sample standard deviation $\sigma_n$ of *CoS* for increasing sample size, *n*, and for Gauss(0). A fitted curve obtained using the least-squares method is also displayed; it is expressed as $\sigma_n = 2.99\, n^{-0.81}$. Similar to bias, very close standard deviation curves are obtained for the 23 copulas used to estimate the bias curve.

TABLE I

SAMPLE MEANS AND SAMPLE STANDARD DEVIATIONS OF THE *CoS* FOR THE GAUSSIAN, GUMBEL, AND CLAYTON COPULA IN THE INDEPENDENCE CASE

| *n* | Gauss(0) $\rho_P = 0$ | | Gumbel(1) $\rho_S = 0$ | | Clayton(0) $\rho_S = 0$ | |
|---|---|---|---|---|---|---|
| | $\mu_n$ | $\sigma_n$ | $\mu_n$ | $\sigma_n$ | $\mu_n$ | $\sigma_n$ |
| 100 | 0.28 | 0.08 | 0.28 | 0.08 | 0.28 | 0.08 |
| 500 | 0.08 | 0.02 | 0.08 | 0.03 | 0.08 | 0.02 |
| 1000 | 0.04 | 0.01 | 0.04 | 0.01 | 0.05 | 0.01 |
| 2000 | 0.02 | 0.01 | 0.02 | 0.01 | 0.02 | 0.01 |
| 3000 | 0.02 | 0.01 | 0.02 | 0.01 | 0.02 | 0.01 |

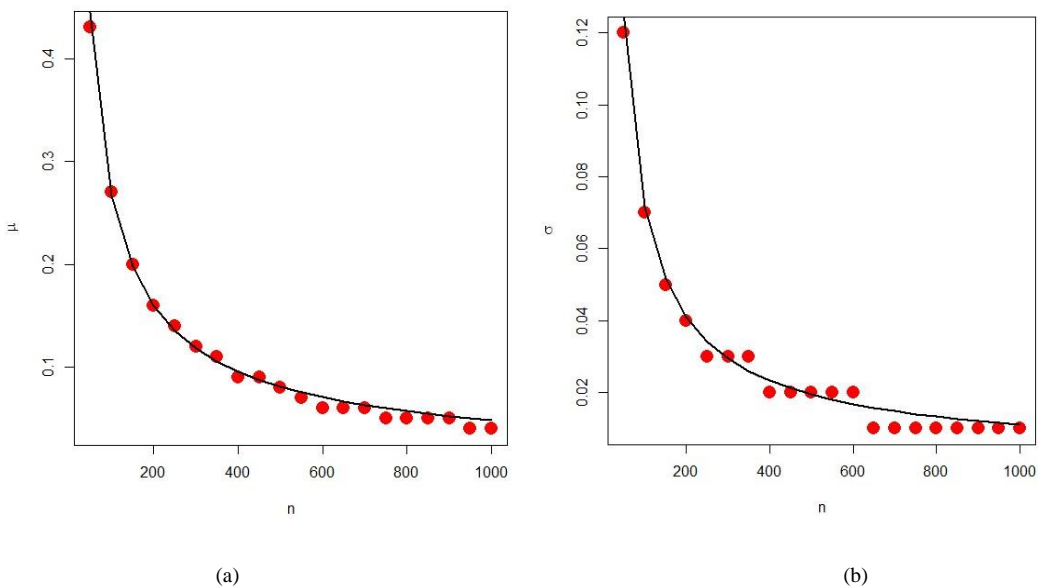

(a)  (b)

**Fig. 7.** (a) Bias mean values and (b) standard deviation values (red solid circles) for the *CoS* along with fitted curves (solid lines) using the least-squares method for the independence case.

## VI. TESTING BIVARIATE INDEPENDENCE

One common practical problem is to test the independence of random variables. To this end, we apply hypothesis testing to the *CoS* based on Corollary 3b). Our intent is to test the null hypothesis, $\mathcal{H}_0$: the random variables are independent, against its alternative, $\mathcal{H}_1$. We standardize the *CoS* under $\mathcal{H}_0$ to get

$$z_n = \frac{CoS - \mu_{n0}}{\sigma_{n0}}, \quad (19)$$

where $\mu_{n0}$ and $\sigma_{n0}$ are the sample mean and the sample standard deviation of the *CoS*, respectively. Resorting to the central limit theorem, we infer that under $\mathcal{H}_0$, $z_n$ approximately follows a standard normal distribution, $\mathcal{N}(0,1)$, for large $n$. This is supported by extensive Monte Carlo simulations that we conduct, where data samples are drawn from various copulas.

As an illustrative example, Fig. 8 displays the QQ-plots of $z_n$ calculated from 100 data sets following Gauss(0.8). It is observed from Fig. 8(a) that the distribution of $z_n$ is skewed to the left for $n = 100$ and from Fig. 8(b) that it is nearly Gaussian for $n = 600$. Hypothesis testing consists of choosing a threshold $c$ at a significance level $\alpha$ under $\mathcal{H}_0$ and then applying the following decision rule: if $|z_n| \leq c$, accept $\mathcal{H}_0$; otherwise, accept $\mathcal{H}_1$. The values of $\mu_{n0}$ and $\sigma_{n0}$ are given by the curves displayed in Fig. 7(a) and 7(b), respectively. Table II displays Type-II errors of the statistical test applied to the *CoS* for Gauss(0) for sample sizes ranging from 100 to 3000. In the simulations, $c = 2.57$ at a significant level of 1% and the alternatives assume weak dependence with $\rho_n = 0.1$ and 0.3. It is observed that Type II-errors decrease as $\rho_n$ increases for a given $n$ and sharply decrease with increasing $n$. This property is related to the Cramer-Rao Lower Bound for correlation, which states that the variance of the estimate decreases as the correlation between the random variables increases [13].

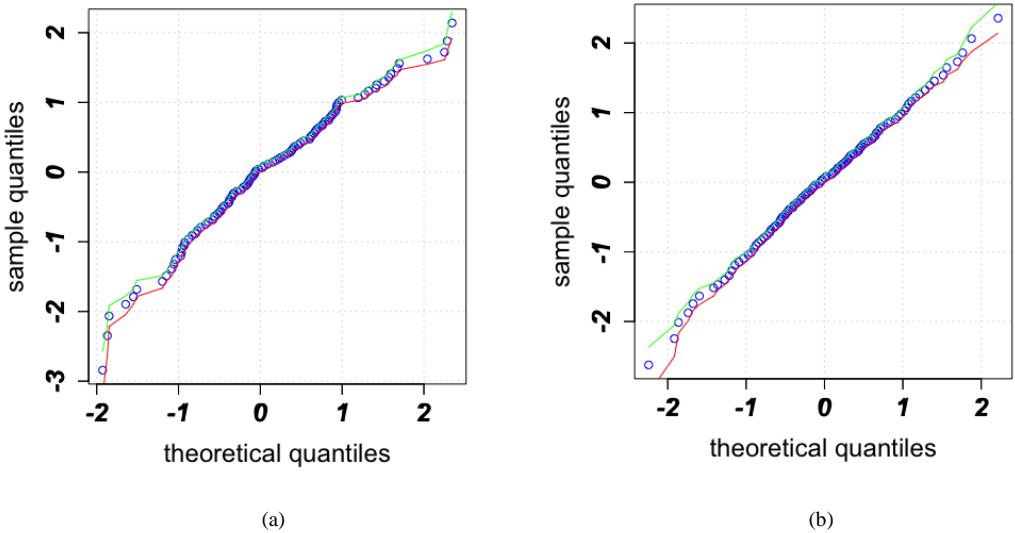

(a)   (b)

**Fig. 8.** Q-Q plots of the standardized *CoS*, $z_n$, based on 100 samples drawn from Gauss(0.8) of size (a) $n = 100$ and (b) $n = 600$. Sample medians and interquartile ranges are displayed in circles and dotted lines, respectively.

TABLE II

TYPE-II ERRORS OF THE STATISTICAL TEST OF BIVARIATE INDEPENDENCE BASED ON *CoS* FOR GAUSS(0)

| $n$ | $\mu_{n0}$ | $\sigma_{n0}$ | Type-II error for $\rho_n = 0.1$ | Type-II error for $\rho_n = 0.3$ |
|---|---|---|---|---|
| 100 | 0.28 | 0.08 | 97% | 46% |
| 500 | 0.08 | 0.02 | 27% | 0% |
| 1000 | 0.04 | 0.01 | 0% | 0% |
| 2000 | 0.02 | 0.01 | 0% | 0% |
| 3000 | 0.02 | 0.01 | 0% | 0% |

## VII. FUNCTIONAL BIVARIATE DEPENDENCE

For monotonic dependence, simulation results show that $CoS = 1$ for all $n \geq 2$. For non-monotonic dependence, there is a bias that becomes negligible when the sample size is sufficiently large. As an illustrative example, Table III displays the sample mean, $\mu_n$, and the sample standard deviation, $\sigma_n$, of the *CoS* for increasing sample size, $n$, for the sinusoidal dependence, $Y = sin(a X)$. It is observed that as the frequency of the sine function increases, the sample bias, $1 - \mu_n$, increases for constant $n$.

TABLE III

SAMPLE MEANS AND SAMPLE STANDARD DEVIATIONS OF THE COS FOR THREE SINUSOIDAL FUNCTIONS OF INCREASING FREQUENCY

| $N$ | $Sin(x)$ | | $Sin(5x)$ | | $Sin(14x)$ | |
|---|---|---|---|---|---|---|
|  | $\mu_n$ | $\sigma_n$ | $\mu_n$ | $\sigma_n$ | $\mu_n$ | $\sigma_n$ |
| 100 | 1.00 | 0.00 | 0.91 | 0.10 | 0.67 | 0.10 |
| 500 | 1.00 | 0.00 | 0.99 | 0.03 | 0.88 | 0.07 |
| 1000 | 1.00 | 0.00 | 1.00 | 0.01 | 0.96 | 0.04 |
| 2000 | 1.00 | 0.00 | 1.00 | 0.00 | 1.00 | 0.01 |
| 3000 | 1.00 | 0.00 | 1.00 | 0.00 | 1.00 | 0.01 |
| 5000 | 1.00 | 0.00 | 1.00 | 0.00 | 1.00 | 0.00 |

## VIII. COPULA-INDUCED DEPENDENCE

Table IV displays $\mu_n$ and $\sigma_n$ of the *CoS* calculated for increasing $n$ and for different degrees of dependencies of two dependent random variables following the Gaussian, the Gumbel and the Clayton copulas. It is interesting to note that for $n \geq 1000$, the *CoS* is nearly equal to the Pearson's $\rho_P$ for the Gaussian copula and to the Spearman's $\rho_S$ for the Gumbel and Clayton copula.

TABLE IV

SAMPLE MEANS AND SAMPLE STANDARD DEVIATIONS OF THE *CoS* FOR THE NORMAL, GUMBEL AND CLAYTON COPULA

| N | Gauss(0.1) $\rho_P = 0.1$ | | Gauss(0.3) $\rho_P = 0.3$ | | Gumbel(1.08) $\rho_S = 0.1$ | | Gumbel(1.26) $\rho_S = 0.3$ | | Clayton(0.15) $\rho_S = 0.1$ | | Clayton(0.51) $\rho_S = 0.3$ | |
|---|---|---|---|---|---|---|---|---|---|---|---|---|
| | $\mu_n$ | $\sigma_n$ | $\mu_n$ | $\sigma_n$ | $\mu_n$ | $\sigma_n$ | $\mu_n$ | $\sigma_n$ | $\mu_n$ | $\sigma_n$ | $\mu_n$ | $\sigma_n$ |
| 100 | 0.33 | 0.09 | 0.49 | 0.09 | 0.34 | 0.09 | 0.51 | 0.10 | 0.34 | 0.09 | 0.49 | 0.09 |
| 500 | 0.14 | 0.05 | 0.36 | 0.05 | 0.16 | 0.05 | 0.37 | 0.05 | 0.15 | 0.05 | 0.37 | 0.05 |
| 1000 | 0.11 | 0.03 | 0.33 | 0.04 | 0.12 | 0.04 | 0.35 | 0.04 | 0.12 | 0.03 | 0.34 | 0.03 |
| 2000 | 0.09 | 0.02 | 0.32 | 0.03 | 0.11 | 0.03 | 0.33 | 0.03 | 0.11 | 0.03 | 0.33 | 0.02 |

### B. Equitability Analysis of the CoS

Reshef *et al.* [7, 8] define the equitability of statistic as its ability to assign equal scores to a collection of functional relationships of the form $(X, Y = f(X))$ subject to the same level of additive noise, which is determined by the coefficient of determination, $R^2$. Here, the noise may be added either to the independent variable $X$, yielding $Y = f(X + \varepsilon)$, or to the response variable, yielding $Y = f(X) + \varepsilon$, or to both, yielding $Y = f(X + \varepsilon_1) + \varepsilon_2$. As defined more formally by Kinney and Atwal [15], a dependence measure $D(X,Y)$ is $R^2$-equitable if and only if, when evaluated on a joint probability distribution $H(x,y)$ that corresponds to a noisy functional relationship between two real random variables $X$ and $Y$, the relation $D(X,Y) = g(R^2[f(X),Y])$. Here, $g(.)$ is a function that does not depend on $H(x,y)$ and $f(.)$ is the function defining the noisy relationship. Interestingly, the authors [15] stressed that no nontrivial measure of dependence can satisfy the mathematical formulation of equitability given above.

We carry out Monte Carlo simulations to assess the $R^2$-equitability of the *CoS* for ten functional relationships of the form $(X, Y = f(X) + \varepsilon)$, which are specified in Table V and where X is drawn from a uniform distribution over [0,1]. In contrast to Reshef *et al.* [32], the noise $\varepsilon$ follows $\mathcal{N}(0, \sigma^2)$, where $\sigma^2 = Var(f(X)) (1/R^2 - 1)$. This expression follows from the definition of $R^2$. Note that $\sigma^2$ is inversely proportional to $R^2$, implying that it increases as $R^2$ approaches zero, that is, as the variability of $Y$ is less and less determined by the variability of $X$. Fig. 9 displays the equitability results of the *CoS* for sample sizes $n = 250, 500$ and $2000$. The red line in Fig. 9 shows the worst interpretable interval, which can be informally defined as the widest range of $R^2$ values corresponding to any one value of the statistic, in this case *CoS*. We observe that the worst and the average interpretable intervals are smaller than one and that they decrease as n increases, indicating an improvement of the equitability of the *CoS* with the sample size. The reader is referred to [32] for a formal definition of the interpretable intervals. We note here that the equitability results depicted in Fig. 9 cannot be compared to those obtained by Reshef *et al.* [7,8] for the *MIC* since the latter are computed under uniform homoscedastic noise on a given interval [15]. We choose to simulate the equitability under Gaussian noise because in signal processing and in communications, Johnson noise follows a Gaussian distribution by virtue of the central limit theorem [33].

TABLE V

FUNCTIONS USED IN THE EQUITABILITY ANALYSIS FOR THE COS

| Function | Color |
|---|---|
| $y = x$ | light blue |
| $y = 4x^2$ | red |
| $y = 41(4x^3 + x^2 - 4x)$ | green |
| $y = \sin(16\pi x)$ | light green |
| $y = \cos(14\pi x)$ | black |
| $y = \sin(10\pi x) + x$ | yellow |
| $y = \sin(6\pi x(1 + x))$ | pink |
| $y = \sin(5\pi x(1 + x))$ | grey |
| $y = 2^x$ | blue |
| $y = 1/10 \sin(10.6 (2x-1)) + 11/10 (2x-1)$ | magenta |

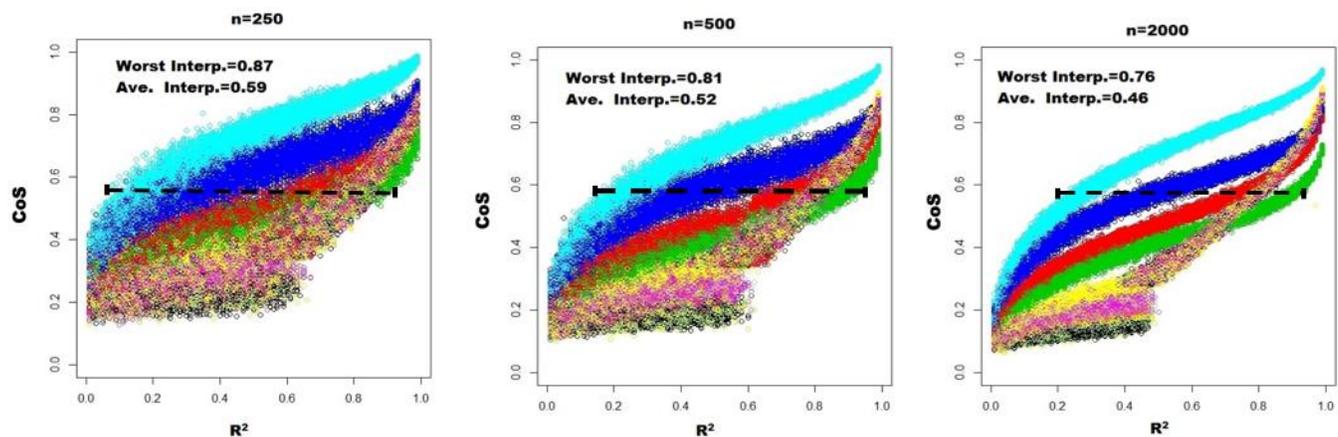

**Fig. 9.** Scatter plots of the *CoS* versus the coefficient of determination, $R^2$, for the ten functional relationships indicated along with their respective colors in Table V and for three sample sizes, n = 250, 500 and 2000. The equitability is measured by means of the worst interpretable interval and the average interpretable interval. The worst interpretable is shown by the dashed red line in the plots above.

## IX. COMPARATIVE STUDY

In this section, we analyze bivariate synthetic datasets and multivariate datasets of real-time stock market returns and of gene regulatory networks.

### A. Bivariate Dependence of Synthetic Data

Let us compare the performances of the *CoS*, *dCor*, *RDC*, *Ccor*, and of the *MICe* for various types of statistical dependencies. Székely *et al.* [12] define the distance correlation, *dCor*, between two random vectors, *X* and *Y*, with finite first moments as

$$dCor(X,Y) = \begin{cases} \frac{v^2(X,Y)}{\sqrt{v^2(X)v^2(Y)}} & \text{for } v^2(X)v^2(Y) > 0, \\ 0 & \text{for } v^2(X)v^2(Y) = 0, \end{cases} \quad (20)$$

where $v^2(X,Y)$ is the distance covariance defined as the norm in the weighted $L_2$ space of the difference between the joint characteristic function and the product of the marginal characteristic functions of *X* and *Y*. Here, $v^2(X)$ stands for $v^2(X,X)$. Lopes-Paz *et al.* [9] define the *RDC* as the largest canonical correlation between *k* randomly chosen nonlinear projections of the copula transformed data. Recall that the largest canonical correlation of two random vectors, ***X*** and ***Y***, is the maximum value of the Pearson correlation coefficient between $a^T X$ and $b^T Y$ for all non-zero real-valued vectors, ***a*** and ***b***. The random nonlinear projections chosen by Lopes-Paz *et al.* [9] are sinusoidal projections with frequencies drawn from the Gaussian distribution. Ding *et al.* [11] define the copula correlation (*Ccor*) as half of the $L_1$ distance between the copula density and the independence copula density. As for the *MIC*, it is defined by Reshef *et al.* [7] as the maximum taken over all *x*-by-*y* grids *G* up to a given grid resolution, typically $x y < n^{0.6}$, of the empirical standardized mutual information, $I_G(x,y)/\log(\min\{x,y\})$, based on the empirical probability distribution over the boxes of a grid G. Formally, we have

$$MIC(X,Y) = \max \left\{ \frac{I_G(x,y)}{\log(\min\{x,y\})} \right\}. \quad (21)$$

In our simulations, we use the *MICe* estimator of the *MIC* when computing a measure of dependence because there is no known algorithm to compute the latter directly in polynomial-time [34]. We also use the *TICe* estimator of the *TIC* when comparing the statistical power. The former is derived by Reshef *et al.* [16] from the *MICe* to achieve high statistical power. Note that we do not include the *TDC* in our analysis because it assumes a given predefined list of functional dependence between the random variables and therefore, its performance strongly depends on the validity of that list; in other words, it is not an agnostic measure of dependence.

*1) Bias analysis for non-functional dependence*

A bias analysis is performed for the *MICe*, the *Ccor*, the *CoS*, the *RDC*, and the *dCor* and using three data samples drawn from a bivariate Gaussian copula with $\rho_P(X,Y) = 0.2$, 0.5 and 0.8, which models a weak, medium and strong dependence, respectively.

The sample sizes range from 50 to 2000, in steps of 50. We observe from Fig. 10 that unlike the *MICe* and *Ccor*, the *CoS, RDC,* and the *dCor* are almost equal to $\rho_P$ for large sample size.

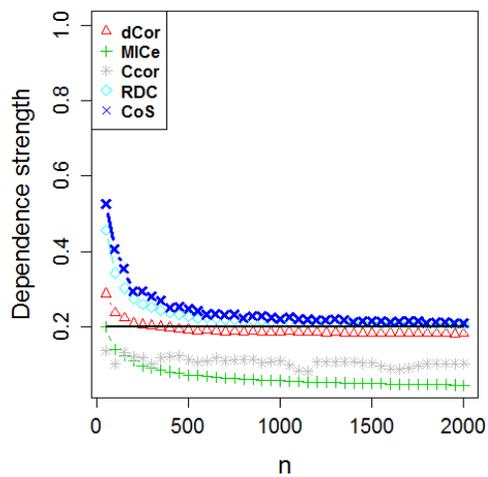

(a)

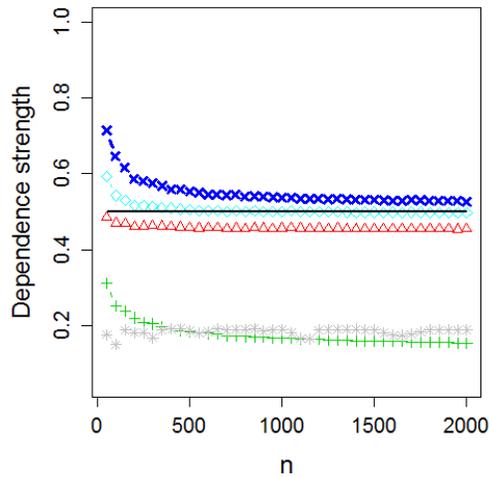

(b)

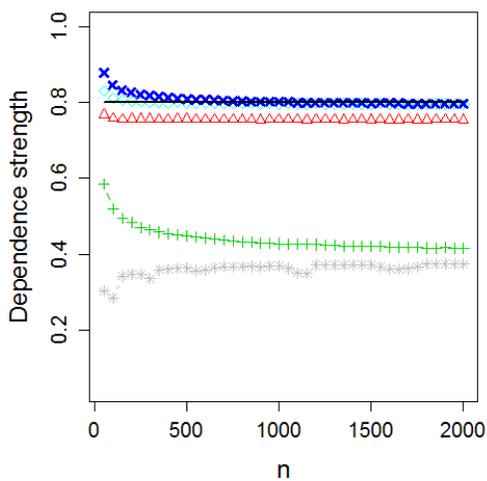

(c)

**Fig. 10**. Bias curves of the *CoS*, *MICe, dCor, RDC,* and *Ccor* for the bivariate Gaussian copula with $\rho_P(X,Y)$ =0.2, 0.5 and 0.8, which are displayed in a), b), and c), respectively, and for sample sizes that vary from 50 and 2000 with steps of 50.

*2) Functional and circular dependence*

In addition to the equitability study reported in Section V.5, we conduct another series of simulations to compare the performance of the *MICe*, the *Ccor*, the *CoS*, the *RDC*, and the *dCor* when they are applied to four data sets drawn from an affine, polynomial, periodic, and circular bivariate relationship with an increasing level of white Gaussian noise. Described in Table VI, the procedure is executed with $N = n = 1000$, where n is the number of realizations of a uniform random variable X and N is the number of times the procedure is executed. We infer from Table VII that while the *CoS*, *dCor*, *Ccor* steadily decrease as the noise level $p$ increases, the *MICe* sharply decreases as $p$ grows from 0.5 to 2 and then reaches a plateau for $p > 2$. This is particularly true for the circular dependence. The *RDC* also decreases steadily with an increase in noise level for the functional dependencies considered, except for the quadratic dependence where it maintains a high power even under heavy noise level.

*3) Ripley's forms and copula's induced dependence*

Table VIII reports values of the *MICe*, the *Ccor*, the *CoS*, the *RDC*, and the *dCor* for Ripley's forms, and copula-induced dependencies for a sample size $n = 1000$ averaged over 1000 Monte-Carlo simulations. The values of the Spearman's $\rho_s$ for Gumbel(5), Clayton(-0.88), Galambos(2), and BB6(2, 2) copulas are calculated using the copula and the CDVine toolboxes of the software package R. As for the four Ripley's forms displayed in Fig. 11, a linear congruential generator using the Box-Muller transformation is used to generate several bivariate sequences with nonlinear dependencies.

TABLE VI
NOISE INSERTION PROCEDURE

*Step 1:* Generate n random realizations of a uniform random variable X on
[-5, 5] to get the data sample $\{x_1, \ldots, x_n\}$;

*Step 2:* Calculate $y_{0i} = f(x_i)$, $i = 1,\ldots,n$, to get n realizations of $Y_0$;

*Step 3:* Replace $y_{0i}$ by $y_{pi}$ for $i = 1,\ldots, n$, according to $y_{pi} = y_{oi}(1 + p\, \varepsilon_i)$,
where $p \in [0, 4]$ and $\varepsilon \sim \mathcal{N}(0, 1)$;

*Step 4*: Calculate *CoS(X,Y)*;

*Step 5*: Repeat Steps 1 through 4 N times and calculate the *CoS* sample mean.

TABLE VII

SAMPLE MEANS OF THE *CoS*, *dCor* AND THE *MICe* FOR SEVERAL DEPENDENCE TYPES AND ADDITIVE NOISE LEVELS

| Type of dependence | Noise level p | 0.5 | 1 | 2 | 3 | 4 |
|---|---|---|---|---|---|---|
| Affine: $Y = 2X+1$ | *CoS* | 0.86 | 0.72 | 0.41 | 0.29 | 0.24 |
| | *dCor* | 0.91 | 0.71 | 0.46 | 0.35 | 0.30 |
| | *MICe* | 0.88 | 0.46 | 0.26 | 0.22 | 0.21 |
| | *RDC* | 0.95 | 0.74 | 0.60 | 0.59 | 0.59 |
| | *Ccor* | 0.63 | 0.47 | 0.34 | 0.30 | 0.29 |
| 4$^{th}$-order Polynomial: $Y=(X^2-0.25)(X^2-1)$ | *CoS* | 0.64 | 0.41 | 0.29 | 0.26 | 0.25 |
| | *dCor* | 0.41 | 0.35 | 0.31 | 0.30 | 0.30 |
| | *MICe* | 0.79 | 0.54 | 0.49 | 0.48 | 0.48 |
| | *RDC* | 0.95 | 0.93 | 0.92 | 0.91 | 0.91 |
| | *Ccor* | 0.72 | 0.63 | 0.60 | 0.59 | 0.59 |
| Periodic: $Y = \cos(X)$ | *CoS* | 0.53 | 0.46 | 0.28 | 0.23 | 0.21 |
| | *dCor* | 0.35 | 0.27 | 0.17 | 0.13 | 0.11 |
| | *MICe* | 0.78 | 0.40 | 0.22 | 0.19 | 0.18 |
| | *RDC* | 0.85 | 0.67 | 0.43 | 0.36 | 0.34 |
| | *Ccor* | 0.57 | 0.41 | 0.29 | 0.26 | 0.24 |
| Circular: $X^2 + Y^2 = 1$ | *CoS* | 0.38 | 0.29 | 0.26 | 0.26 | 0.26 |
| | *dCor* | 0.12 | 0.10 | 0.10 | 0.10 | 0.10 |
| | *MICe* | 0.13 | 0.09 | 0.08 | 0.08 | 0.08 |
| | *RDC* | 0.51 | 0.38 | 0.36 | 0.36 | 0.36 |
| | *Ccor* | 0.20 | 0.15 | 0.14 | 0.14 | 0.14 |

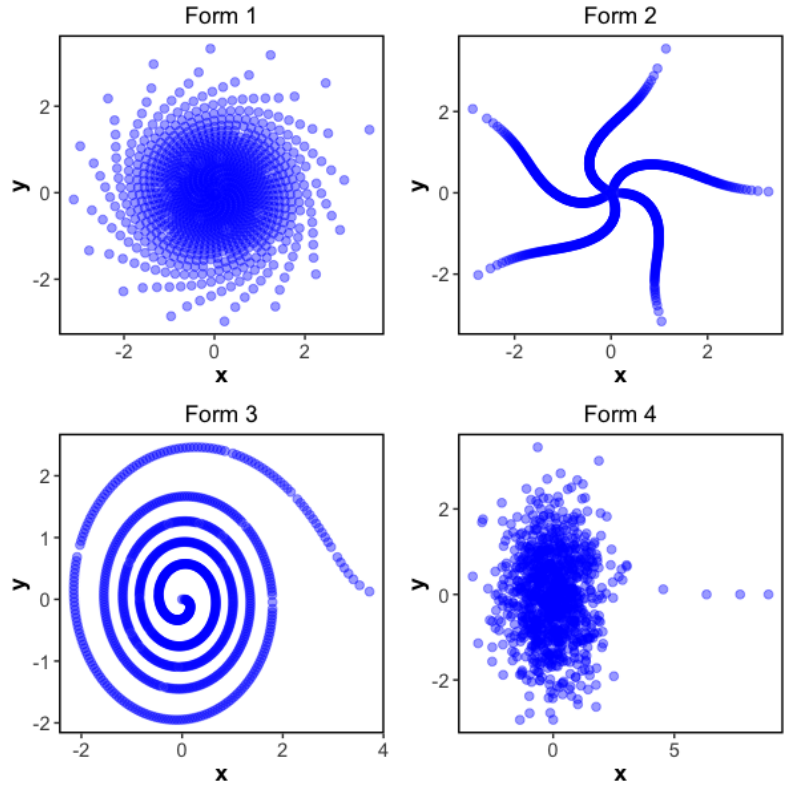

**Fig. 11.** Four Ripley's plots generated using a linear congruential generator followed by the Box-Muller transformation. The parameters of the congruential generator, $x_{i+1} = (a\, x_i + c)$ modulo $M$, are as follows: Form 1: $a = 65$, $c = 1$, $M = 2048$; Form 2: $a = 1229$, $c = 1$, $M = 2048$; Form 3: $a = 5$, $c = 1$, $M = 2048$; Form 4: $a = 129$, $c = 1$, $M = 2^{64}$.

Table VIII shows that the *CoS*, *MICe*, *RDC*, and *Ccor* correctly reveal some degree of nonlinear dependence for Ripley's form 2, with the *Ccor* detecting the highest level of dependence and the *dCor* the lowest level. It is observed that the *Ccor* is the only metric to correctly reveal some degree of nonlinear dependence for Ripley's form 3. Furthermore, unlike the *MICe* values, the *dCor* and the *CoS* values are very close to the Pearson's $\rho_P$ value for the Gaussian copula and to the Spearman's $\rho_S$ values for the Gumbel, Clayton, Galambos and BB6 copulas. This does not come as a surprise because, as proved in [35], $\rho_P(X,Y)$ and $\rho_S(X,Y)$ can be respectively expressed in terms of the copula of $X$ and $Y$ as

$$\rho_P(X,Y) = \frac{1}{\sigma_X \sigma_Y} \iint_0^1 [C(u,v) - uv] dF_1^{-1}(u) dF_2^{-1}(v), \qquad (22)$$

and

$$\rho_S(X,Y) = 12 \iint_0^1 [C(u,v) - uv] du dv. \qquad (23)$$

Here $\sigma_X$ and $\sigma_Y$ denote the standard deviation of $X$ and $Y$, respectively. Noting the similarity between these relationships and the expression of the distance function, $\lambda(C(u,v))$, given by Definition 5, we conjecture that asymptotically, $CoS(X,Y) = \rho_P(X,Y)$ for the Gaussian copula and $CoS(X,Y) = \rho_S(X,Y)$ for the other above-mentioned copulas.

TABLE VIII

DEPENDENCE INDICES FOR COPULA-INDUCED DEPENDENCIES AND RIPLEY'S FORMS FOR A SAMPLE SIZE $n = 1000$

| Type of Dependence | CoS | dCor | MICe | RDC | Ccor |
|---|---|---|---|---|---|
| Ripley's form 1 | 0.01 | 0.02 | 0.02 | 0.02 | 0.01 |
| Ripley's form 2 | 0.52 | 0.19 | 0.42 | 0.42 | 0.84 |
| Ripley's form 3 | 0.14 | 0.08 | 0.12 | 0.13 | 0.26 |
| Ripley's form 4 | 0.03 | 0.04 | 0.03 | 0.08 | 0.09 |
| Gaussian(0.1), $\rho_P = 0.10$ | 0.11 | 0.10 | 0.04 | 0.13 | 0.10 |
| Gumbel(5), $\rho_S = 0.94$ | 0.92 | 0.93 | 0.72 | 0.96 | 0.62 |
| Clayton(-0.88), $\rho_S = -0.87$ | 0.90 | 0.87 | 0.68 | 0.88 | 0.75 |
| Galambos(2), $\rho_S = 0.81$ | 0.82 | 0.79 | 0.48 | 0.86 | 0.42 |
| BB6(2,2), $\rho_S = 0.80$ | 0.84 | 0.83 | 0.57 | 0.92 | 0.48 |

[1] $Y = -X - 1$ for $-5 \leq X \leq -1$; $Y = X + 1$ for $-1 \leq X \leq 0$; $Y = -X + 1$ for $0 \leq X \leq 1$; $Y = X - 1$ for $1 \leq X \leq 5$.

*B. Statistical Power Analysis*

Finally, following Simon and Tibshirani [36], we investigate the power of the statistical tests based on the *CoS*, *dCor*, *RDC*, *TICe*, and the *Ccor* for bivariate independence subject to increasing additive Gaussian noise levels. Recall that the power of a test is defined as the probability to accept the alternative hypothesis $\mathcal{H}_1$ when $\mathcal{H}_1$ is true. The procedure implemented is described in Table IX. Six noisy functional dependencies at a noise level $p$ ranging from 10% to 300% are considered. They include a linear, a quadratic, a cubic, a fourth-root, a sinusoidal, and a circular dependence. Fig. 12 displays the power of the tests calculated using a collection of $N = 500$ data samples, each of size $n = 500$, for a significance level $\alpha = 5\%$ under the null hypothesis. As observed from that figure, the *CoS* is a powerful measure of dependence for the linear, cubic, circular and rational dependence, but it is less powerful for the quadratic and periodic cases. More specifically, the *CoS* performs better than any other copula-based metric for the linear and fourth root dependencies, performs equally well to the *RDC* for the cubic and circular dependencies, but performs less than the *RDC* for the quadratic and sinusoidal dependencies. We conjecture that for the latter dependencies the loss of power of the *CoS* is due to the decrease of performance of the procedure for finding the local optima described in Step 7 of Section IV.C as the noise level increases. As compared to the non-copula based metrics, the *CoS* performs equally well or better in the linear, cubic, circular, and rational dependence cases, but performs slightly worse than the *dCor* in the quadratic and sinusoidal dependencies. The performance of the *CoS* estimator and improvements to this algorithm for quadratic and sinusoidal dependencies is an area of future research.

TABLE IX

PROCEDURE FOR COMPUTING THE POWER OF A STATISTICAL TEST OF BIVARIATE INDEPENDENCE

*Step 1:* Under the null hypothesis, $\mathcal{H}_0$: $X$ and $Y = f(X_1)$ are independent, generate $N$ samples of size $n$ of two independent random variables, $X \sim U(0,1)$ and $X_1 \sim U(0,1)$, and calculate $Y = f(X_1) + p\,\varepsilon$, where $\varepsilon \sim \mathcal{N}(0,1)$ and $p \in [0,3]$;

*Step 2:* Calculate the $N$ values of a statistic (*e.g.*, the *CoS, TICe, RDC, dCor, Ccor*) and pick the 95 percentile as the cutoff value, $c$;

*Step 3:* Under the alternative, $\mathcal{H}_1$: $X$ and $Y = f(X)$ are dependent, generate $N$ samples of size $n$ of $X \sim U(0,1)$ and calculate $Y = f(X) + p\,\varepsilon$, where $\varepsilon \sim \mathcal{N}(0,1)$;

*Step 4*: Calculate the $N$ values of the statistic;

*Step 5*: Calculate the power $P_d$ of the test given by

$$P_d = \frac{1}{N}\sum_{i=1}^{N} \mathbf{1}\{CoS_i > c,\ i = 1, \dots, N\}$$

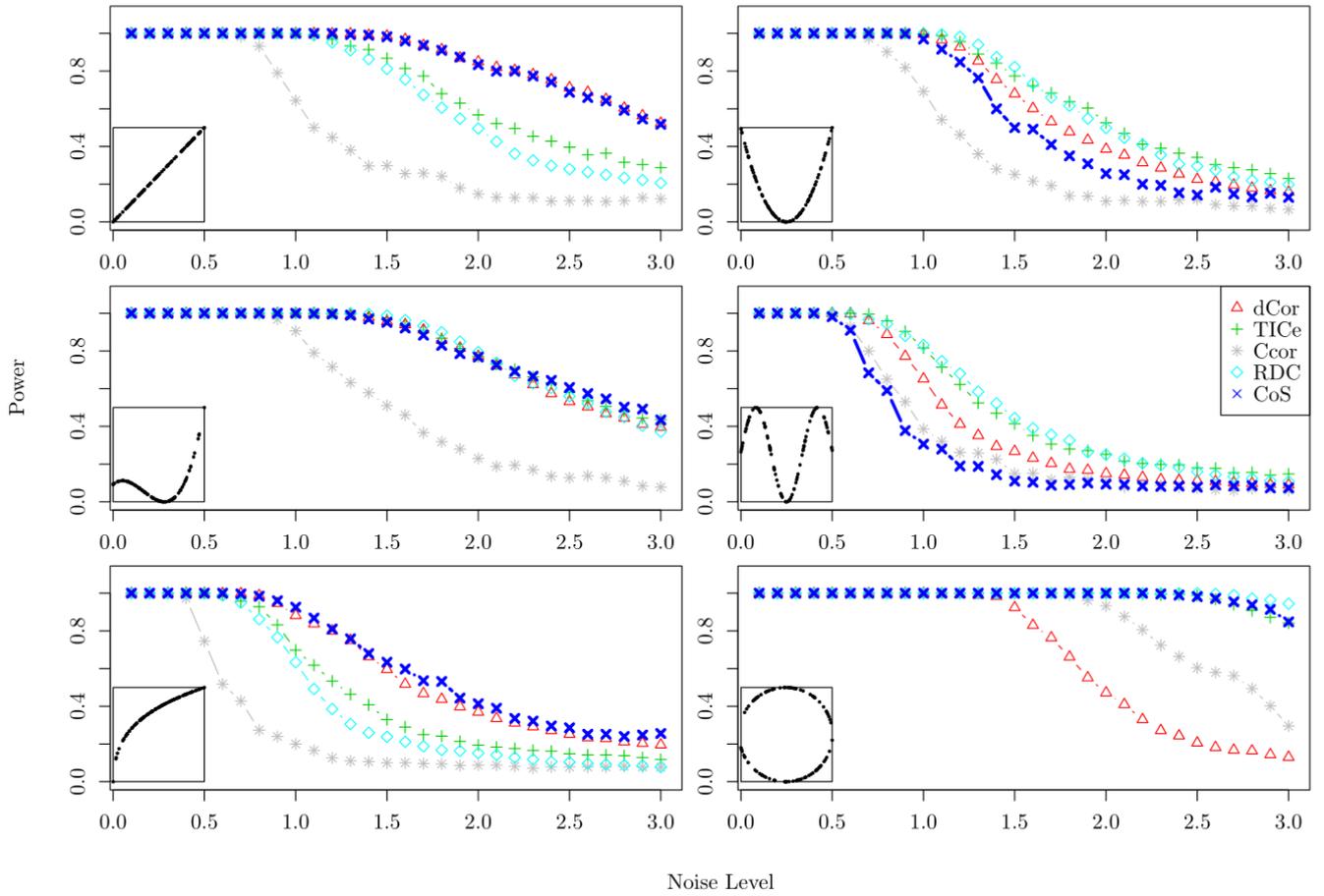

**Fig. 12**. Power of five statistical tests of bivariate independence based on the *CoS* (blue cross), the *dCor* (red triangle), the *Ccor* (grey star), the *RDC* (teal diamond), and the *TICe* (green plus) calculated from $N = 500$ data samples of size $n = 500$ and a significance level of $\alpha = 5\%$ under the null hypothesis for six noisy functional dependencies with a noise level $p$ ranging from 0.0 to 3.0. The dependency type for each plot is shown in the bottom left inlet plot for each subplot.

## C. Dependence Analyses of Data Sets from Realistic Systems

### 1) Dependence Analysis of Real-Time Stock Market Index Returns

There exists a large literature in finance [37-39], especially in stock market [40, 41], that deals with nonlinear dependence between stochastic signals. Here, we analyze the dependence between the returns of four stock market indices recorded monthly from January 1991 to November 2016, namely the Standard and Poor's (S&P) 500 index, the Deutscher AktienindeX (DAX) index, the Nikkei 225, and the CAC 40 index. The price data was acquired from Yahoo Finance, and returns were calculated by computing a one-lag time difference of the price data [41]. Some bivariate scatter plots of these indices returns are displayed in Fig. 13. As observed, the two strongest dependencies are between S&P 500 and CAC 40 with an approximate linear dependence, and between DAX and CAC 40 with an approximate linear dependence.

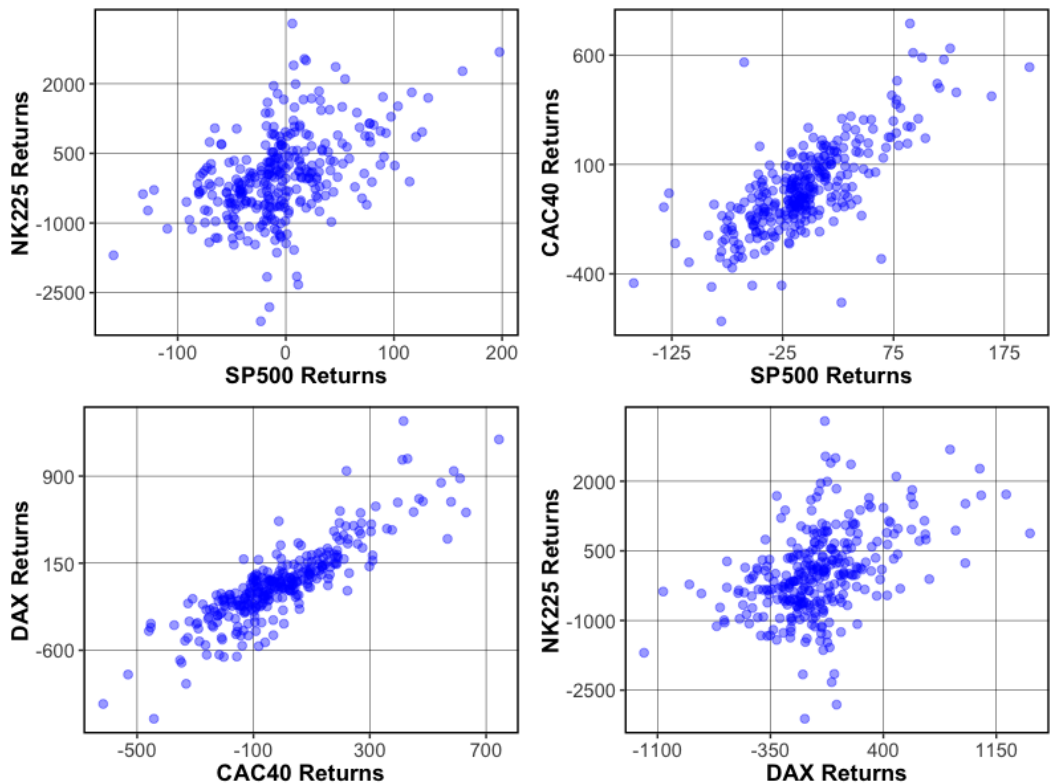

**Fig. 13**. Scatter plots of four stock markets index returns, namely Nikkei 225 vs. S&P 500; DAX vs. CAC 40; CAC 40 vs. S&P 500; Nikkei 225 vs. S&P 500.

Before computing the pairwise and multivariate dependencies between the four stock market returns described above, we first test the stationarity of these time-series signals using the augmented Dickey-Fuller test [42]. This technique tests the null hypothesis that a unit-root is present in a time-series sample, which is an indication that the time-series is not stationary. For the four time-series signals used in the analysis, the test yields a p-value of less than 0.01. Even at a significance level of $\alpha=99\%$, the null hypothesis is rejected and the alternative hypothesis that the series is stationary is accepted. We use this result as a basis of our analysis for measuring the dependence between the stock market returns data.

Table X provides values taken by the *CoS, MICe, RDC,* and the *Ccor* for the six pairwise dependencies between these stock market returns. We do not include the *dCor* in this comparison because it assumes independent and identically distributed samples [12], whereas the stock market returns analyzed here are typically serially correlated. While all the statistics show evidence of some degree of dependence between the six pairs of index returns, they greatly differ in their values for the two aforementioned strongest dependent cases. Indeed, for the S&P 500-Nikkei 225 pair, the *CoS* and the *RDC* both exhibit values close to 0.5, while the *MICe* takes a low value of 0.2. Considering the scatter plot in Fig. 13, the *CoS* and the *RDC* provide more a credible view of dependence between the returns of these indices. As for the S&P 500 and CAC 40 pair, the *CoS* and the *RDC* take values close to 0.8 while the *MICe* attains 0.38. From the scatter plots in Fig. 13, the values reached by the *CoS*, and the *RDC* seem more reasonable than the value taken by the *MICe*. Additionally, the *Ccor* consistently underestimates all pairwise dependencies for the returns data.

Now, the question that arises is to know whether a pairwise dependence analysis of these time series is sufficient to unveil all the dependencies between them. The answer is negative as shown in Table XI. Indeed, the trivariate dependence analysis using the *CoS* reveals a strong dependence between S&P 500, CAC 40 and DAX with a value of 0.72 and a weaker dependence between the other triplets. The multivariate dependence assessment between the four index returns, namely S&P500, CAC 40, Nikkei 225, and DAX indicates a medium dependence with a *CoS* value of 0.55. Recall that such analysis cannot be carried out with the *MICe,* the *dCor,* or the *RDC* since they are restricted to bivariate dependence.

TABLE X

COMPARISON OF BIVARIATE DEPENDENCE VALUES BETWEEN FOUR STOCK MARKET INDEX RETURNS

| Pairwise Returns | CAC40 | | | | DAX | | | | Nikkei 225 | | | |
|---|---|---|---|---|---|---|---|---|---|---|---|---|
| | *CoS* | *MICe* | *RDC* | *Ccor* | *CoS* | *MICe* | *RDC* | *Ccor* | *CoS* | *MICe* | *RDC* | *Ccor* |
| **S&P 500** | 0.76 | 0.38 | 0.82 | 0.30 | 0.75 | 0.37 | 0.82 | 0.34 | 0.46 | 0.20 | 0.51 | 0.14 |
| **CAC 40** | - | - | - | - | 0.85 | 0.61 | 0.89 | 0.43 | 0.43 | 0.19 | 0.49 | 0.12 |
| **DAX** | - | - | - | - | - | - | - | - | 0.52 | 0.20 | 0.52 | 0.20 |

TABLE XI

*CoS* Values Of Multivariate Dependence Between Three and Four Stock Market Index Returns

| Multivariate Dependence Between Three and Four Stock Market Index Returns | CoS |
|---|---|
| S&P500, Nikkei 225, and DAX | 0.57 |
| S&P500, CAC40, and DAX | 0.75 |
| S&P500, CAC40, and NIKKEI225 | 0.56 |
| CAC40, Nikkei 225, and DAX | 0.58 |
| S&P500, CAC40, Nikkei 225, and DAX | 0.58 |

*2) Dependence Analysis of Gene Expression Data*

Gene regulatory networks are other examples where nonlinear dependence between gene expression signals is prominent [43-45]. Large databases of these signals obtained from microarray technique assays [46] are posted in the Internet for processing, with the aim to construct genome control maps that reveal how regulatory genes (i.e., hubs) activates or repress regulated genes.

Following the procedure implemented in the Minet package developed by Meyer *et al.* [47], we calculate the levels of pairwise dependence of the gene expression signals given by the *CoS, dCor, RDC, Ccor* and the *MICe*, which form the elements of the dependence matrix, **M**, except for the diagonal elements, which are set to zero. Note that we assume here that a gene expression is a continuous random variable as recommended in [48]. Then, we take in turn each element of the matrix **M** as the threshold of a binary statistical test that is applied to all the elements of **M**. The decision rule is a follows: there is no link between a pair of genes if the *CoS* or the *MICe* value is less than the threshold; otherwise, there is a link. Next, by comparing the inferred links between pairs of genes and the supposedly true links of the gene regulatory network, we estimate the True Positive (TP), the False Positive (FP), the True Negative (TN), and the False Negative (FN) defined as the number of positive (respectively negative) decisions over the total number of decisions. Finally, we calculate the maximum of the F-scores calculated over all the elements of the matrix **M** and plot the Receiver Operating Characteristic (ROC) curves. Recall that the F-score is defined as

$$F-score = \frac{2\ Precision\ x\ Recall}{Precision\ +Recall}, \qquad (24)$$

where Precision = TP/(TP + FP) and Recall = TP/(TP + FN).

The foregoing procedure is applied to two sets of gene expression signals. The first set consists of synthetic signals of fifty genes of the yeast genome collected from 100 experiments utilizing the microarray data generator SynTReN [49]. It has been extracted from the Internet databases [47] by executing the Minet package [47]. The ROC curves depicted in Fig. 14(a) and 14(b) show that the *CoS* exhibits similar performance to the other metrics compared. This is confirmed by the ROC areas and the F-score maximum values shown in Table XII.

The procedure is also applied to a second set of real signals of eight genes of the *E. Coli* SOS response pathway to DNA damage, whose expressions are assumed to occur according to the true regulatory network displayed in Fig. 15 [50]. Retrieved from [51], the set consists of 196 data per gene, which include small proportions of missing values. The ROC curves displayed in Fig. 16(a) and 16(b) show that the *CoS* performs similarly to the other metrics, except for the *Ccor*, which performs poorly for this dataset. This is in agreement with the ROC areas and the F-score maximum values shown in Table XII.

TABLE XII

ROC AREAS AND MAXIMA OF THE F-SCORES FOR THE *MICe* AND *CoS* FOR YEAST AND E. COLI GENE EXPRESSION DATA SETS

| Performance Indices | | *MICe* | *CoS* | *RDC* | *dCor* | *Ccor* |
|---|---|---|---|---|---|---|
| Yeast Genes | ROC area | 0.79 | 0.72 | 0.83 | 0.81 | 0.81 |
| | F-score-max | 0.40 | 0.42 | 0.46 | 0.34 | 0.44 |
| E. Coli Genes | ROC area | 0.85 | 0.74 | 0.82 | 0.86 | 0.32 |
| | F-score-max | 0.63 | 0.55 | 0.67 | 0.74 | 0.37 |

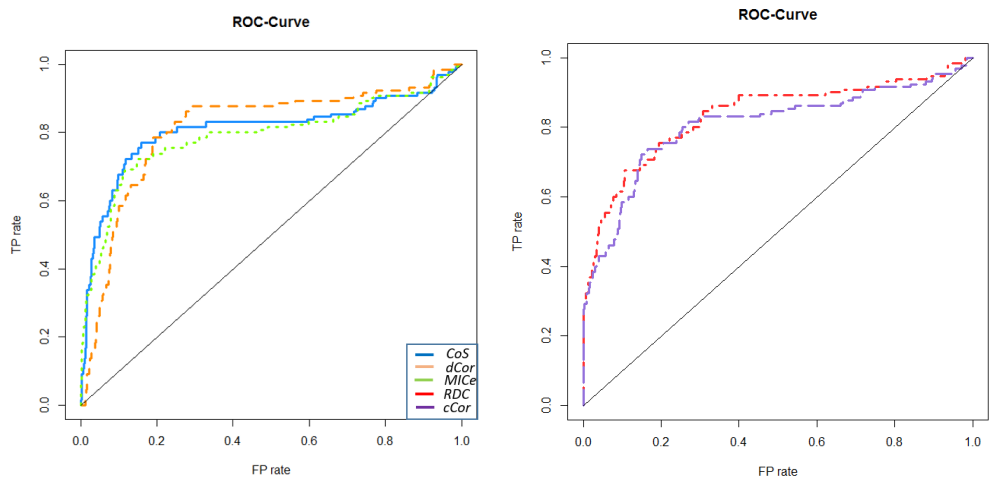

**Fig. 14**. ROC curves of the *CoS* and the *MICe* of 50 gene expression signals of the yeast genome.

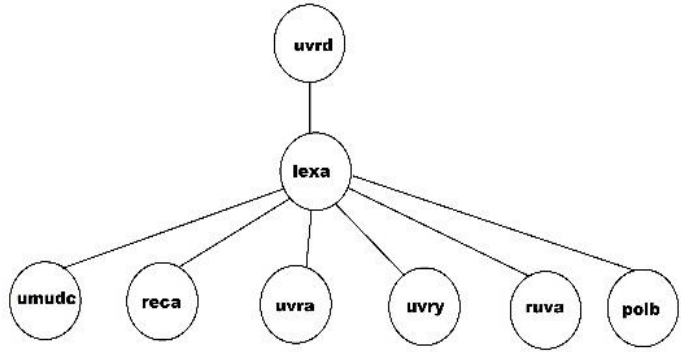

**Fig. 15.** True regulatory subnetwork of eight *E. Coli* genes for protein production. Lexa is a regulatory gene that activates or represses seven genes.

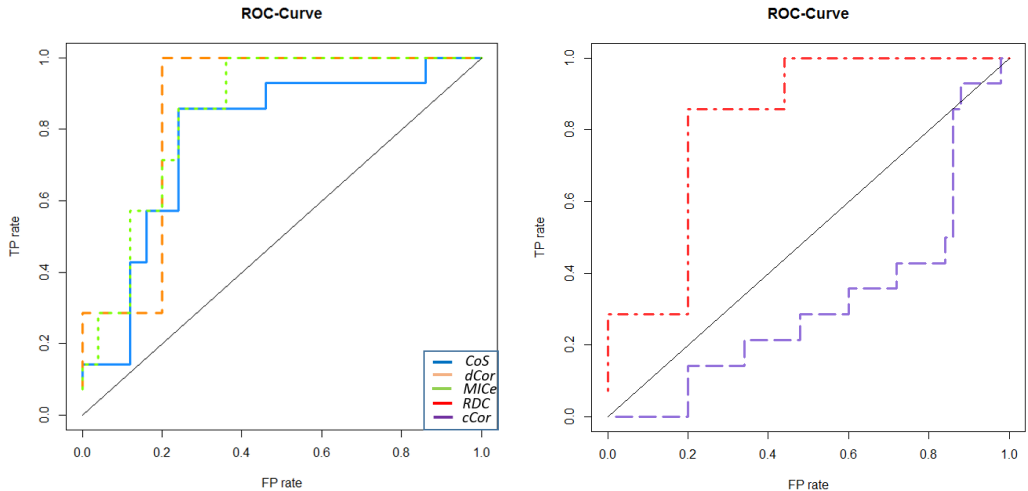

**Fig. 16**. ROC curves of the *CoS* and the *MICe* of eight gene expressions of the *E. Coli* SOS response pathway to DNA damage.

VII. CONCLUSIONS AND FUTURE WORK

A new statistic for multivariate nonlinear dependence, the *CoS*, has been proposed and its statistical properties unveiled. In particular, it asymptotically approaches zero for statistical independence and one for functional dependence. Finite-sample bias and standard deviation curves of the *CoS* have been estimated and hypothesis testing rules have been developed to test bivariate independence. The power of the *CoS*-based test and its $R^2$-equitability has been evaluated for noisy functional dependencies. Monte Carlo simulations show that the *CoS* performs reasonably well for both functional and non-functional dependence and exhibits a good power for testing independence against all alternatives. By virtue of Theorem 2.6 proved in Embrechts *et al.* [38], it follows that the *CoS* is invariant to strictly increasing functional transforms; other invariance properties of the *CoS* will be investigated as a future work. Another interesting property of the *CoS* that is not shared by the *MICe, RDC*, *Ccor*, and the *dCor* is its ability to measure multivariate dependence. This has been demonstrated using stock market index returns. Good performance of the *CoS* has been shown in gene expressions of regulatory networks. Note that the code that implements the *CoS* is available on the GitHub repository [52]. As a future research work, we will assess the self-equitability of the *CoS* and other metrics under various noise probability distributions, including thick tailed distributions such as the Laplacian distribution and long memory processes, and we will investigate the robustness of the *CoS* to outliers. Furthermore, we will apply the *CoS* to common signal processing and machine learning problems, including data mining, cluster analysis, and testing of independence.

ACKNOWLEDGEMENTS

The authors are grateful to David N. Reshef for sending them the Java package that implements the *TICe*.

**Mohsen Ben Hassine** received an engineering diploma and the M.S. degree in computer sciences from the École Nationale des Sciences de l' Informatique, Tunis, Tunisia, in 1993 and 1996, respectively. He is currently a graduate teaching assistant at the University of El Manar, Tunis, Tunisia. His research interests include statistical signal processing, mathematical simulations, and statistical bioinformatics.

**Lamine Mili** received an electrical engineering diploma from the Swiss Federal Institute of Technology, Lausanne, in 1976, and a Ph.D. degree from the University of Liege, Belgium, in 1987. He is presently a Professor of Electrical and Computer Engineering at Virginia Tech. His research interests include robust statistics, robust statistical signal processing, radar systems, and power system analysis and control. Dr. Mili is a Fellow of IEEE for contribution to robust state estimation for power systems.

**Kiran Karra** received a B.S. in Electrical and Computer Engineering and an M.S. degree in Electrical Engineering from North Carolina State University and Virginia Polytehcnic Institute and State University, in 2007 and 2012, respectively. He is currently a research associate at the Virginia Tech and is studying statistical signal processing and machine learning for his PhD research.